\documentclass[11pt]{article}
\usepackage[a4paper, total={6in, 8in}]{geometry}
\usepackage{enumerate}
\usepackage{amsmath}
\usepackage{amssymb,latexsym}
\usepackage{amsthm}
\usepackage{color}
\usepackage{cancel}
\usepackage{graphicx}
\usepackage{mathtools}

\usepackage{hyperref}
\usepackage{comment}
\usepackage{amscd}
\usepackage{cite}
\usepackage{lineno}

\makeatletter
\usepackage{comment}
\let\wfs@comment@comment\comment
\let\comment\@undefined

\usepackage{changes}
\let\wfs@changes@comment\comment
\let\comment\@undefined

\newcommand\comment{%
    \ifthenelse{\equal{\@currenvir}{comment}}
    {\wfs@comment@comment}
    {\wfs@changes@comment}%
}

\definechangesauthor[name=Daniele, color=red]{DAN}
\definechangesauthor[name=Giovanni, color=blue]{GIO}
\definechangesauthor[name=Nico, color=green]{NICO}

\newtheorem{theorem}{Theorem}[section]

\newtheorem{definition}[theorem]{Definition}

\newtheorem{open}[theorem]{Open Problem}
\newtheorem*{theorem*}{Theorem}

\newcommand\PG{\mathrm{PG}}

\newcommand\mO{\mathcal{O}}

\title{Ovoids of $Q^+(7,q)$ of low-degree}
\author{Daniele Bartoli\thanks{Department of Mathematics and Informatics, University of Perugia, Perugia, Italy.
Email address: {\tt{daniele.bartoli@unipg.it}}},\, Nicola Durante\thanks{Department of Mathematics and Applications R. Caccioppoli, University of Naples Federico II, Naples, Italy.
Email address: {\tt{ndurante@unina.it}}},\, Giovanni Giuseppe Grimaldi\thanks{Department of Mathematics and Applications R. Caccioppoli, University of Naples Federico II, Naples, Italy.
Email address: {\tt{giovannigiuseppe.grimaldi@unina.it}}},\, Marco Timpanella \thanks{Department of Mathematics and Informatics, University of Perugia, Perugia, Italy.
Email address: {\tt{marco.timpanella@unipg.it}}}}
\date{}
\begin{document}

\maketitle

\begin{abstract}
Ovoids of the hyperbolic quadric $Q^+(7,q)$ of $\mathrm{PG}(7, q)$ have been extensively studied over the past 40 years, partly due to their connections with other combinatorial objects. It is well known that the points of an ovoid of $Q^+(7,q)$ can be parametrized by three polynomials $f_1(X,Y, Z)$, $f_2(X,Y, Z)$, $f_3(X,Y, Z)$. In this paper, we classify ovoids of $Q^+(7,q)$ of low degree, specifically under the assumption that $f_1(X,Y, Z)$, $f_2(X,Y, Z)$, $f_3(X,Y, Z)$ have degree at most $3$.
Our approach relies on the analysis of an algebraic hypersurface associated with the ovoid.
\end{abstract}
\section{Introduction}\label{Sec:Intro}
Let $\mathbb{F}_{q^n}$, where $q$ is a power of a prime $p$, be the finite field with $q^n$ elements. A \textit{finite classical polar space}, denoted as $\mathbb{P}$, is the set of absolute points associated with either a polarity or a non-degenerate quadratic form within the projective space $\mathrm{PG}(m, q)$. The largest projective subspaces fully contained in $\mathbb{P}$ are referred to as its \textit{generators}.

\begin{definition}
\label{ovoids}
An \textit{ovoid} of a finite classical polar space $\mathbb{P}$ in $\mathrm{PG}(m, q)$ is a subset of points in $\mathbb{P}$ such that every generator of $\mathbb{P}$ intersects the ovoid in exactly one point.
\end{definition}

In this work, we focus on the study of ovoids of the hyperbolic quadric $Q^+(7, q)$. Consider the projective space $\mathrm{PG}(7, q)$ with homogeneous coordinates $(X_0, X_1, X_2, X_3, X_4, X_5, X_6, X_7)$. Let $H_\infty$ represent the hyperplane at infinity, defined by $X_0 = 0$. 

There are exactly two types of non-degenerate quadrics in $\mathrm{PG}(7, q)$ up to collineations: the elliptic quadric $Q^-(7, q)$ and the hyperbolic quadric $Q^+(7, q)$. Here, we focus on the hyperbolic quadric defined by the equation:
\[
Q: X_0X_7 + X_1X_6 + X_2X_5 + X_3X_4 = 0.
\]

An ovoid of $Q^+(7, q)$ consists of $q^3 + 1$ mutually non-collinear points lying entirely on the quadric. Without loss of generality, one can assume that an ovoid $\mathcal{O}_7$ contains the points $(1, 0, 0, 0, 0, 0, 0, 0)$ and $P_\infty = (0, 0, 0, 0, 0, 0, 0, 1)$. Consequently, any ovoid of $Q^+(7, q)$ can be parameterized as:
\begin{align}
\label{Eq:Param}
\mathcal{O}_7 = &\big\{(1, x, y, z, f_1(x, y, z), f_2(x, y, z), f_3(x, y, z),\\
&\hspace*{2 cm}-z f_1(x, y, z) - y f_2(x, y, z) - x f_3(x, y, z)) \, \mid \, x, y, z \in \mathbb{F}_q\big\} \cup \{P_\infty\}\nonumber,
\end{align}
where $f_1, f_2, f_3 : \mathbb{F}_q^3 \to \mathbb{F}_q$ are functions such that $f_i(0, 0, 0) = 0$ for $i \in \{1, 2, 3\}$. This ovoid is denoted by $\mathcal{O}_7(f_1, f_2, f_3)$, and its \textit{degree} is defined as $d = \max\{\deg(f_1), \deg(f_2), \deg(f_3)\}$.

For simplicity, a point of $\mathcal{O}_7(f_1, f_2, f_3) \setminus \{P_\infty\}$ is denoted by $P(x, y, z)$. The set $\mathcal{O}_7(f_1, f_2, f_3)$ forms an ovoid if and only if
\begin{equation}
\label{Param}
\langle P(x_1, y_1, z_1), P(x_2, y_2, z_2) \rangle \neq 0,
\end{equation}
for all distinct $(x_1, y_1, z_1)$ and $(x_2, y_2, z_2)$ in $\mathbb{F}_q^3$, where $\langle \cdot, \cdot \rangle$ is the bilinear form associated with the quadratic form defining $Q^+(7, q)$. This bilinear form is symmetric if $q$ is odd, and alternating if $q$ is even.

Since $\langle P(x, y, z), (0, 0, 0, 0, 0, 0, 0, 1) \rangle = 1$, Condition \eqref{Param} can be rewritten as:
\begin{align}
\label{eqov}
&(x_1 - x_2)(f_3(x_2, y_2, z_2) - f_3(x_1, y_1, z_1)) \nonumber \\
&+ (y_1 - y_2)(f_2(x_2, y_2, z_2) - f_2(x_1, y_1, z_1)) \nonumber \\
&+ (z_1 - z_2)(f_1(x_2, y_2, z_2) - f_1(x_1, y_1, z_1)) \neq 0,
\end{align}
for all $(x_1, y_1, z_1) \neq (x_2, y_2, z_2)$ in $\mathbb{F}_q^3$.

Previous research \cite{OvoidiQ+5, OvoidiQ6, OvoidiQ4} has employed techniques from algebraic geometry to classify ovoids of quadrics in lower dimensions, such as $Q(4, q)$, $Q^+(5, q)$, and $Q(6, q)$. For the case of $Q^+(7, q)$, the known ovoids are summarized in Table \ref{Table1}, with more details available in \cite{Williams}.

It is known that $Q^+(7,q)$ has a unique ovoid if $q \in \{2,3,4\}$, see \cite{Gun, Jho, Pat}.

This paper aims to partially classify ovoids of $Q^+(7, q)$ by examining specific algebraic varieties associated with the ovoid. For a survey of the interplay between algebraic varieties and combinatorial objects, see \cite{BartoliSurvey}. Here, we restrict our study to ovoids parameterized by \eqref{Eq:Param}, where each $f_i$ is of degree at most 3. 

The paper is organized as follows.In Section \ref{Sec:AlgebraicVarieties}, we explore the relationship between the ovoids of $Q^+(7,q)$ and certain  algebraic hypersurfaces of $\mathrm{PG}(6,q)$, introducing the algebraic geometry tools necessary to study this connection.

In Section \ref{sec:known} we list the triples of functions \( f_1, f_2, f_3 \) associated with some known ovoids, following the description in Equation \eqref{Eq:Param}. In particular, we provide the functions \( f_1, f_2, f_3 \) that describe the Kantor ovoids for $q=2^h$, $h\geq 1$ and the Kantor ovoid for $q=p^h$, $p\equiv 3 \pmod{3}$, $h$ odd; see Section \ref{Kantorq=2} and \ref{subsec:kantor2}. As a consequence, we show that these two families of ovoids are parameterized by polynomials of degree $2$ and $3$, respectively.

Finally, in Sections \ref{sec:classdeg2} and \ref{classdeg3}, we address the classification problem for ovoids of degrees $2$ and $3$. More specifically, in Section \ref{sec:classdeg2}, we prove that if $q$ is sufficiently large, the only ovoid of degree $2$ is the Kantor ovoid for $q=2^h$, $h\geq 1$; see Theorem \ref{Th:d=2}. In the case of degree $3$, we provide a complete classification of the functions \( f_1, f_2, f_3 \) that describe ovoids of $Q^+(7,q)$, assuming that $q$ is sufficiently large and $p\neq 3$; see Theorems \ref{Th:d=3,quadriche1}, \ref{Th:d=3,quadriche2}, and \ref{Th:d=3,quadriche3}. 

\begin{center}
\begin{table}[h!]
\caption{Known ovoids of $Q^+(7, q)$}\label{Table1}
\begin{center}
\begin{tabular}{||c c||}
\hline
Name & Restrictions \\ [0.5ex] 
\hline\hline
\textit{Thas-Kantor} & $q = 3^h$, $h > 0$ \\ 
\hline
\textit{Ree-Tits} & $q = 3^{2h+1}$, $h > 0$ \\  
\hline
\textit{Kantor} & $q = p^h$, $p \equiv 2 \pmod{3}$ (prime), $h$ odd \\ 
\hline
\textit{Kantor} & $q = 2^h$, $h \geq 1$ \\
\hline
\textit{Dye} & $q = 8$ \\  
\hline
\textit{Conway et al./Moorhouse} & $q \geq 5$ and prime \\
\hline
\end{tabular}
\end{center}
\end{table}
\end{center}

\section{Link between ovoids of $Q^+(7,q)$ and algebraic varieties}\label{Sec:AlgebraicVarieties}
 This section aims to explore the relationship between the ovoids of $Q^+(7,q)$ and specific algebraic hypersurfaces, introducing the mathematical tools necessary to study this connection.

An algebraic hypersurface $\mathcal{S}$ refers to an algebraic variety described by a single equation. If the polynomial defining a hypersurface over a field $\mathbb{K}$ cannot be factored over any algebraic extension of $\mathbb{K}$, it is called \emph{absolutely irreducible}. For a hypersurface $\mathcal{S}$ defined by a polynomial $F$, an absolutely irreducible $\mathbb{K}$-rational component is an irreducible hypersurface whose polynomial has coefficients in $\mathbb{K}$ and divides $F$. Throughout this work, we will rely on the use of homogeneous equations for algebraic varieties. For those interested in a more thorough introduction to algebraic varieties, we recommend \cite{Hartshorne}.

To derive results on non-existence and partial classifications for the ovoids $\mathcal{O}_7(f_1,f_2,f_3)$ of $Q$, we focus on the hypersurface $\mathcal{S}_{f_1,f_2,f_3} \subset \mathrm{PG}(6,q)$ given by the homogeneous equation $F(X_0,X_1,X_2,X_3,X_4,X_5,X_6) = 0$, where
\begin{eqnarray}\label{Eq:F}
  F(1,X_1,X_2,X_3,X_4,X_5,X_6) &=& (X_1 - X_4) \Big( f_3(X_4,X_5,X_6) - f_3(X_1,X_2,X_3) \Big) \nonumber \\
  && + (X_2 - X_5) \Big( f_2(X_4,X_5,X_6) - f_2(X_1,X_2,X_3) \Big) + \nonumber \\
  && (X_3 - X_6) \Big( f_1(X_4,X_5,X_6) - f_1(X_1,X_2,X_3) \Big).
\end{eqnarray}

According to Condition \eqref{Param}, $\mO_7(f_1,f_2,f_3)$ is an ovoid of $Q$ if and only if $\mathcal{S}_{f_1,f_2,f_3}$ does not contain any affine $\mathbb{F}_q$-rational points outside the 3-space defined by $X_1 - X_4 = X_2 - X_5 = X_3 - X_6 = 0$.

A crucial result for determining the presence of rational points in algebraic varieties over finite fields is the Lang-Weil theorem \cite{LangWeil} (1954), which generalizes the Hasse-Weil bound.
\begin{theorem}\label{Th:LangWeil}[Lang-Weil Theorem]
Let $\mathcal{V} \subset \PG(n,q)$ be an absolutely irreducible variety of dimension $r$ and degree $d$. There exists a constant $C$, depending only on $n$, $r$, and $d$, such that the following inequality holds:
\begin{equation}\label{Eq:LW}
\left|\#(\mathcal{V} \cap \mathrm{PG}(n,q)) - \sum_{i=0}^{r} q^i \right| \leq (d-1)(d-2)q^{r-1/2} + Cq^{r-1}.
\end{equation}
\end{theorem}

\noindent While the constant $C$ is not explicitly given in \cite{LangWeil}, numerous explicit estimates have been derived in works such as \cite{CafureMatera, Ghorpade_Lachaud, Ghorpade_Lachaud2, LN1983, WSchmidt, Bombieri}, often of the form $C = f(d)$ where $q > g(r,d)$ and $f$ and $g$ are low-degree polynomials. For a survey of these bounds, see \cite{CafureMatera}. Further excellent surveys on results of the Hasse-Weil and Lang-Weil type are provided in \cite{Geer, Ghorpade_Lachaud2}. A specific result by Cafure and Matera \cite{CafureMatera} is given below.

\begin{theorem}\cite[Theorem 7.1]{CafureMatera}\label{Th:CafureMatera}
Let $\mathcal{V} \subset \mathrm{AG}(n,q)$ be an absolutely irreducible $\mathbb{F}_q$-variety of dimension $r > 0$ and degree $d$. If $q > 2(r+1)d^2$, the following estimate holds:
$$ |\#(\mathcal{V} \cap \mathrm{AG}(n,q)) - q^r| \leq (d-1)(d-2)q^{r-1/2} + 5d^{13/3} q^{r-1}. $$
\end{theorem}

For the purposes of this work, the presence of an absolutely irreducible $\mathbb{F}_q$-rational component in $\mathcal{S}_{f_1,f_2,f_3}$ is sufficient to derive asymptotic results on non-existence.

\begin{theorem}\label{Th:Main}
Let $\mathcal{S}_{f_1,f_2,f_3}: F(X_0,X_1,X_2,X_3,X_4,X_5,X_6) = 0$, where $F$ is defined as in \eqref{Eq:F}. Assume that $q > 6.3 (d+1)^{13/3}$, $\max\{\deg(f_1), \deg(f_2), \deg(f_3)\} = d$, and that $\mathcal{S}_{f_1,f_2,f_3}$ contains an absolutely irreducible component $\mathcal{V}$ defined over $\mathbb{F}_q$. Then, $\mO_7(f_1,f_2,f_3)$ cannot be an ovoid of $Q$.
\end{theorem}
\proof
Given that $q > 6.3 (d+1)^{13/3}$, by Theorem \ref{Th:CafureMatera}, with $r = 5$, it follows that $\mathcal{V}$ contains more than $q^3$ affine $\mathbb{F}_q$-rational points. Thus, at least one affine $\mathbb{F}_q$-rational point $(a_1, b_1, c_1, a_2, b_2, c_2)$ must lie in $\mathcal{V} \subset \mathcal{S}_{f_1,f_2,f_3}$ such that $a_1 \neq a_2$ or $b_1 \neq b_2$ or $c_1 \neq c_2$, which implies that $\mO_7(f_1,f_2,f_3)$ is not an ovoid of $Q$.
\endproof

\section{Known Ovoids and Their Associated Functions}\label{sec:known}

In this section, we list the triples of functions \( f_1, f_2, f_3 \) associated with some known ovoids, following the description in Equation \eqref{Eq:Param}. 

To this end, recall that to any point \( P(x,y,z) \) of an ovoid \( O \) of \( Q \), there corresponds the 3-space \( S_P \) of \( Q \) given by:
\begin{eqnarray}
&& X_0 = x X_6 - y X_5 + z X_4, \nonumber \\
&& X_1 = -x X_7 + f_1(x,y,z) X_5 + f_2(x,y,z) X_4, \nonumber \\
&& X_2 = y X_7 - f_1(x,y,z) X_6 + f_3(x,y,z) X_4, \nonumber \\
&& X_3 = -z X_7 - f_2(x,y,z) X_6 - f_3(x,y,z) X_5, \nonumber
\end{eqnarray}
and to the point \( P_\infty \), there corresponds the 3-space \( X_4 = X_5 = X_6 = X_7 = 0 \).

Then, Condition \eqref{eqov}, which states that \( O \) is an ovoid, is equivalent to the condition that the corresponding 3-spaces form a spread of \( Q \). Moreover, since the ovoid contains the point \( P(0,0,0) \), the corresponding spread contains the 3-space \( X_0 = X_1 = X_2 = X_3 = 0 \).

This implies that among the following matrices

\[
\left\{
\begin{pmatrix}
0 & x & -y & z \\
-x & 0 & f_1(x,y,z) & f_2(x,y,z) \\
y & -f_1(x,y,z) & 0 & f_3(x,y,z) \\
-z & -f_2(x,y,z) & -f_3(x,y,z) & 0
\end{pmatrix} : x,y,z \in \mathbb{F}_q
\right\}
\]

all, except the zero matrix, are non-singular. 

This means that to a spread of \( Q \), there corresponds a so-called Kerdock set of \( q^3 \) skew-symmetric matrices of order 4 such that the difference of any two matrices has full rank. This is also a maximum rank distance code with \( d = 4 \) in the set of skew-symmetric matrices of order 4, see \cite{DelGoe}.

\subsection{Kantor Ovoid \( q = 2^h \), \( h \geq 1 \)}\label{subsec:kantor2}

In \cite[page 90]{Williams}, it is proved that for \( q \in \{2,4,16\} \), the Kantor ovoid is described by the functions
\begin{eqnarray}
    && f_1(x,y,z) = xy + z^2, \nonumber \\
    && f_2(x,y,z) = xz + y^2 + z^2, \nonumber \\
    && f_3(x,y,z) = yz + x^2 + y^2 + z^2. \nonumber
\end{eqnarray}
Moreover, it is observed that when \( q = 8 \), these functions do not define an ovoid.

In this section, we provide a description of the Kantor ovoid in the general case \( q = 2^h \), with \( h \geq 1 \).

Let \( \mathbb{V} = \mathbb{F}_q \oplus \mathbb{F}_{q^3} \oplus \mathbb{F}_{q^3} \oplus \mathbb{F}_q \) equipped with 
\[
Q((a, \gamma, \delta, d)) = ad + \mathrm{Tr}_{q^3/q}(\gamma \delta).
\]

Then, the Kantor ovoid is given by the set
\begin{eqnarray}\label{KO}
\{(0,0,0,1)\} \cup \{(1,t,t^{q+q^2}, \mathrm{N}_{q^3/q}(t)): t \in \mathbb{F}_{q^3}\},
\end{eqnarray}
see e.g. \cite[page 82]{Williams}. Let \( \{1, \alpha, \beta\} \) be an \( \mathbb{F}_q \)-basis of \( \mathbb{F}_{q^3} \), and put \( t = x + y \alpha + z \beta \) with \( x, y, z \in \mathbb{F}_q \). Hence,
\[
t^{q+q^2} 
= x^2 + (\alpha^q + \alpha^{q^2}) xy + (\beta^q + \beta^{q^2}) xz + 
y^2 \alpha^{q+q^2} + yz(\alpha^{q^2} \beta^q + \alpha^q \beta^{q^2}) + z^2 \beta^{q+q^2}.
\]

Now, by direct computations, it is possible to check that \( t^{q+q^2} \) equals
\[
x f_3(x,y,z) + y f_2(x,y,z) + z f_1(x,y,z),
\]
where 
\begin{eqnarray}
f_1(x,y,z) &=& \mathrm{Tr}_{q^3/q}(\alpha \beta^2 + \alpha \beta^{q^2}) xy + \mathrm{Tr}_{q^3/q}(\beta) x^2 + \mathrm{Tr}_{q^3/q}(\alpha^{q+1} \beta^{q^2}) y^2 + \mathrm{N}_{q^3/q}(\beta) z^2, \nonumber \\
f_2(x,y,z) &=& \mathrm{Tr}_{q^3/q}(\alpha \beta^2 + \alpha \beta^{q^2}) xz + \mathrm{Tr}_{q^3/q}(\alpha) x^2 + \mathrm{N}_{q^3/q}(\alpha \beta^{q^2}) y^2 + \mathrm{Tr}_{q^3/q}(\alpha \beta^{q^2+q}) z^2, \nonumber \\
f_3(x,y,z) &=& \mathrm{Tr}_{q^3/q}(\alpha \beta^2 + \alpha \beta^{q^2}) yz + x^2 + \mathrm{Tr}_{q^3/q}(\alpha^{q+1} \beta^{q^2}) y^2 + \mathrm{Tr}_{q^3/q}(\beta^{q+1}) z^2. \label{Eq:Kantor_q_pari}
\end{eqnarray}

The ovoid in Equation \eqref{KO} corresponds to the set
\[
\{ (1, x, y, z, f_1(x,y,z), f_2(x,y,z), f_3, z f_1 + y f_2 + x f_3) : x, y, z \in \mathbb{F}_q \} \cup \{(0,0,0,0,0,0,0,1)\}.
\]

\subsection{Thas-Kantor, Ree-Tits, and Dye Ovoids}

Starting with the first two ovoids in Table \ref{Table1}, we recall that these are ovoids of the 6-dimensional parabolic quadric \( Q(6,q) \). Thus, we easily derive the functions for the Thas-Kantor ovoid, which are the following:
$$f_1(x,y,z) = z, \quad f_2(x,y,z) = -\mu y^3 + x^2 y - xz, \quad  f_3(x,y,z) = -\frac{1}{\mu} x^3 + x y^2 + yz,$$
where \( \mu \) is a non-square, and the functions for the Ree-Tits ovoid are the following:
$$f_1(x,y,z) = z, f_2(x,y,z) = -x^{\sigma+3} + y^\sigma + x^2 y - xz,  f_3(x,y,z) = -x^{2\sigma+3} + x^\sigma y^\sigma - z^\sigma + x y^2 + yz,$$
where \( \sigma = \sqrt{3q} \).

Finally, the Dye ovoid is defined for \( q = 8 \), and the functions \( f_1, f_2, f_3 \) are the following:
\begin{eqnarray}
&& f_1(x,y,z) = x + y + z + x^2 y + x^4 y^2 + x y^2 + x^2 y^4 + x^4 y^4, \nonumber \\
&& f_2(x,y,z) = y + x^2 z + x^4 z^2 + x z^2 + x^2 z^4 + x^4 z^4, \nonumber \\
&& f_3(x,y,z) = x + y + y^2 z + y^4 z^2 + y z^2 + y^2 z^4 + y^4 z^4, \nonumber
\end{eqnarray}
see \cite[page 89]{Williams}.

\subsection{Kantor ovoid $q\equiv 2 \pmod{3}$}\label{Kantorq=2}

\textbf{Case $q$ odd}\\
Since $q \equiv 2 \pmod 3$, $-3$ is not a square in $\mathbb{F}_q$. Let $\xi \in \mathbb{F}_{q^2} \setminus \mathbb{F}_q$ be such that $\xi^2 = -3$. Consider the vector space 
$$\mathbb{V} = \left\{
\begin{pmatrix}
x & y & c \\
z & a & y^q \\
b & z^q & x^q
\end{pmatrix} : x, y, z \in \mathbb{F}_{q^2}, a, b, c \in \mathbb{F}_q, a + x + x^q = 0
\right\}.$$
Let $x = x_0 + \xi x_1$, $y = y_0 + \xi y_1$, $z = z_0 + \xi z_1$ so that $\mathrm{Tr}_{q^2/q}(x) = 2x_0$, $\mathrm{N}_{q^2/q}(x) = x_0^2 + 3x_1^2$, and $\mathrm{Tr}_{q^2/q}(yz) = 2(y_0 z_0 - 3y_1 z_1)$. Now, $\mathbb{V}$ is isomorphic (as an $\mathbb{F}_q$-vector space) to the 8-dimensional vector space 
$$\mathbb{V}' := \left\{
\left(b, x_0 + x_1, z_0, z_1, -6y_1, 2y_0, 3(x_0 - x_1), c\right) : x_0, x_1, y_0, y_1, z_0, z_1, b, c \in \mathbb{F}_q \right\};$$
see \cite[page 76]{Williams}.

The quadratic form 
$$Q(M) = \mathrm{Tr}_{q^2/q}(x)^2 - \mathrm{N}_{q^2/q}(x) + \mathrm{Tr}_{q^2/q}(yz) + bc$$
reads
$$3x_0^2 - 3x_1^2 + 2(y_0 z_0 - 3y_1 z_1) + bc,$$
which is $X_0 X_7 + X_1 X_6 + X_2 X_5 + X_3 X_4$, in terms of the coordinates of $\mathbb{V}'$.

Thus, the ovoid corresponds to 
\begin{eqnarray*}
    \Omega & := & \{(1, x_0 + x_1, z_0, z_1, -6(x_1 z_0 - x_0 z_1), 2(x_0 z_0 + 3x_1 z_1), 3(x_0 - x_1), x_0^2 + 3x_1^2) \\
    && : x_0, x_1, z_0, z_1 \in \mathbb{F}_q, 2x_0 + z_0^2 + 3z_1^2 = 0 \} \cup \{(0, 0, 0, 0, 0, 0, 0, 1)\},
\end{eqnarray*}
see \cite[page 82]{Williams}. Let $\alpha = x_0 + x_1$, $\beta = z_0$, $\gamma = z_1$. Note that 
$$x_0 = -\frac{z_0^2 + 3z_1^2}{2} = -\frac{\beta^2 + 3\gamma^2}{2}, \qquad x_1 = \alpha - x_0 = \frac{2\alpha + \beta^2 + 3\gamma^2}{2}.$$
Thus,
\begin{eqnarray*}
f_1 &=& -6\left(\frac{2\alpha \beta + \beta^3 + 3\gamma^2 \beta}{2} + \frac{\beta^2 \gamma + 3\gamma^3}{2}\right) \\
    &=& -6\alpha \beta - 3\beta^3 - 9\gamma^3 - 9\gamma^2 \beta - 3\beta^2 \gamma, \\
f_2 &=& 2\left(-\frac{\beta^3 + 3\beta \gamma^2}{2} + 3\frac{2\alpha \gamma + \beta^2 \gamma + 3\gamma^3}{2}\right) \\
    &=& -\beta^3 - 3\beta \gamma^2 + 6\alpha \gamma + 3\beta^2 \gamma + 9\gamma^3, \\
f_3 &=& -3\alpha - 3\beta^2 - 9\gamma^2,
\end{eqnarray*}

and the left-hand side in Condition \eqref{eqov} reads
\begin{eqnarray*}
3\left(U + x_5 W - x_6 V + x_5 V + 3 x_6 W + \frac{-\xi + 3}{6} V^2 + \frac{-\xi + 3}{2} W^2\right) \cdot \\
\left(U + x_5 W - x_6 V + x_5 V + 3 x_6 W + \frac{\xi + 3}{6} V^2 + \frac{\xi + 3}{2} W^2\right),
\end{eqnarray*}
where $U = x_1 - x_4$, $V = x_2 - x_5$, $W = x_3 - x_6$. It is readily seen that the above quantity vanishes only when $U = V = W = 0$. \\

\noindent \textbf{Case $q$ even}\\
Since $q \equiv 2 \pmod{3}$, there is $\xi \in \mathbb{F}_{q^2} \setminus \mathbb{F}_q$ such that $\xi^2 = 1 + \xi$ and $\mathrm{Tr}_{q^2/q}(\xi) = 1$. Using the same notations as in the odd case, we have $\mathrm{Tr}_{q^2/q}(x) = x_1$, $\mathrm{N}_{q^2/q}(x) = x_0^2 + x_0 x_1 + x_1^2$, and $\mathrm{Tr}_{q^2/q}(yz) = y_0 z_1 + y_1 z_0 + y_1 z_1$. Then the vector space $\mathbb{V}$ is isomorphic to the 8-dimensional $\mathbb{F}_q$-vector space 
$$
\mathbb{V}' = \{ (b, x_0 + x_1, z_1, z_0 + z_1, y_1, y_0, x_0, c) : b, c, x_0, x_1, y_0, y_1, z_0, z_1 \in \mathbb{F}_q \}.
$$

The quadratic form 
$$Q(M) = \mathrm{Tr}_{q^2/q}(x)^2 + \mathrm{N}_{q^2/q}(x) + \mathrm{Tr}_{q^2/q}(yz) + bc$$
reads now
$$x_0(x_0 + x_1) + y_0 z_1 + y_1 (z_0 + z_1) + bc,$$
which is $X_0 X_7 + X_1 X_6 + X_2 X_5 + X_3 X_4$,  again with respect to the coordinates of $\mathbb{V}'$. The ovoid corresponds to
\begin{eqnarray*}
    \Omega &=& \{(1, x_0 + x_1, z_1, z_0 + z_1, x_0 z_1 + x_1 z_0, x_0 (z_0 + z_1) + x_1 z_1, x_0, x_0^2 + x_0 x_1 + x_1^2) \\
    && : x_0, x_1, z_0, z_1 \in \mathbb{F}_q, x_1 = z_0^2 + z_0 z_1 + z_1^2 \} \cup \{(0, 0, 0, 0, 0, 0, 0, 1)\}.
\end{eqnarray*}
Let $\alpha = x_0 + x_1$, $\beta = z_1$, $\gamma = z_0 + z_1$. Note that 
$$
z_0 = \gamma + \beta, \quad 
x_1 = \gamma^2 + \gamma \beta + \beta^2, \quad 
x_0 = \alpha + \gamma^2 + \gamma \beta + \beta^2, 
$$
and therefore
\begin{eqnarray*}
f_1 &=& \left( \alpha + \gamma^2 + \gamma \beta + \beta^2 \right) \beta + \left( \gamma^2 + \gamma \beta + \beta^2 \right) \left( \gamma + \beta \right) \\
    &=& \gamma^3 + \gamma^2 \beta + \gamma \beta^2 + \alpha \beta, \\
f_2 &=& \left( \alpha + \gamma^2 + \gamma \beta + \beta^2 \right) \gamma + \left( \gamma^2 + \gamma \beta + \beta^2 \right) \beta \\
    &=& \gamma^3 + \beta^3 + \alpha \gamma, \\
f_3 &=& \alpha + \gamma^2 + \gamma \beta + \beta^2.
\end{eqnarray*}
Thus, the left-hand side in Condition \eqref{eqov} reads
\begin{eqnarray*}
\left( U + x_5 (V + W) + x_6 W + \xi (V^2 + VW + W^2) \right) \cdot \\
\left( U + x_2 (V + W) + x_3 W + \xi (V^2 + VW + W^2) \right),
\end{eqnarray*}
where $U = x_1 + x_4$, $V = x_2 + x_5$, $W = x_3 + x_6$. Since $V^2 + VW + W^2 = (V + \xi W)^{q+1}$,  the above quantity vanishes  if and only if $U = V = W = 0$.

\section{Classification results for ovoids of degree $2$}\label{sec:classdeg2}

The aim of this section is to classify, for $q$ large enough, the ovoids $\mO_7(f_1,f_2,f_3)$ of $Q$ under the assumption that $f_1,f_2,f_3$ are of degree $2$.

To this end, let
\begin{equation}\label{Eq:f_i2}
f_i(X,Y,Z)=A_iXY+B_iXZ+C_i YZ+D_i X^2+E_i Y^2+F_i Z^2,
\end{equation}
with $A_i, B_i, C_i,D_i, E_i,F_i\in \mathbb{F}_q$ and $i\in \{1,2,3\}$.

 By Theorem \ref{Th:Main}, if $\mO_7(f_1,f_2,f_3)$ is an ovoid of $Q$, then the hypersurface $\mathcal{S}_{f_1,f_2,f_3}: F(X_0,X_1,X_2,X_3,X_4,X_5,X_6)=0$ does not contain any absolutely irreducible component defined over $\mathbb{F}_q$, and hence it must split into three hyperplanes switched by the action of the $q$-Frobenius automorphism $\varphi_q$, acting as follows on the polynomials defining the components
 $$ \sum a_{i_1,\ldots,i_6}X_1^{i_1}X_2^{i_2}X_3^{i_3} X_4^{i_4}X_5^{i_5}X_6^{i_6} \mapsto \sum a_{i_1,\ldots,i_6}^qX_1^{i_1}X_2^{i_2}X_3^{i_3} X_4^{i_4}X_5^{i_5}X_6^{i_6}. $$ Moreover, $\mO_7(f_1,f_2,f_3)$ is an ovoid of $Q$ if and only if the affine points of $\mathcal{S}_{f_1,f_2,f_3}$ are contained in  the 3-space $X_1-X_4=X_2-X_5=X_3-X_6=0$, so the three hyperplanes must be defined by $\pi=0$, $\pi_q=0$, and $\pi_{q^2}=0$, where
\begin{eqnarray*}
\pi&=&    (X_1-X_4)+\alpha(X_2-X_5)+\beta(X_3-X_6)=0\\
\pi_q&=&  (X_1-X_4)+\alpha^q(X_2-X_5)+\beta^{q}(X_3-X_6)=0\\
\pi_{q^2}&=&  (X_1-X_4)+\alpha^{q^2}(X_2-X_5)+\beta^{q^2}(X_3-X_6)=0,
\end{eqnarray*}
and $\{1,\alpha,\beta\}$ is a basis of $\mathbb{F}_{q^3}$ over $\mathbb{F}_q$. 

Therefore, $F(1,X_1,X_2,X_3,X_4,X_5,X_6)$ factors as $\pi \cdot \pi_q \cdot \pi_{q^2}$, and $G:=F-\pi \cdot \pi_q \cdot \pi_{q^2}$ is the zero polynomial. In particular, this shows that $q$ is even as the coefficients of $X_1^3$ and $X_4^3$ in $G$ are $D_3-1$ and $D_3+1$, respectively, and they must vanish. Moreover, since all the other coefficients of $G$ vanish (see Appendix for the complete list), we have 
$$\begin{array}{c}
    A_2=A_3=B_1=B_3=C_1=C_2=0,\\
    A_1=B_2=C_3=\mathrm{Tr}_{q^3/q}(\alpha\beta^q + \alpha \beta^{q^2}),\\ 
    D_1=\mathrm{Tr}_{q^3/q}(\beta), \quad    D_2=\mathrm{Tr}_{q^3/q}(\alpha), \quad
    D_3 =1\\
     E_1=\mathrm{Tr}_{q^3/q}(\alpha^{q+1}\beta^{q^2}), \quad
     E_2=\alpha^{q^2+q+1}, \quad
     E_3=\mathrm{Tr}_{q^3/q}(\alpha^{q+1})\\
      F_1=\beta^{q^2+q+1}, \quad
    F_2=\mathrm{Tr}_{q^3/q}(\alpha\beta^{q^2+q}), \quad
     F_3=\mathrm{Tr}_{q^3/q}(\beta^{q+1}).\\
\end{array}$$

Observe that with this choice, all the coefficients $A_i, B_i, C_i, D_i, E_i, F_i$ are in $\mathbb{F}_q$, as requested, and $f_1,f_2,f_3 $ are as in \eqref{Eq:Kantor_q_pari}. Therefore, we have the following.
\begin{theorem}\label{Th:d=2}
Let $q>6.3\cdot 3^{13/3}$. If
\begin{align*}
\label{Eq:Param}
\mO_7(f_1,f_2,f_3):=&\big\{(1, x, y, z, f_1(x, y, z), f_2(x, y, z), f_3(x, y, z),\\
&\hspace*{1.5 cm}-z f_1(x, y, z) - y f_2(x, y, z) - x f_3(x, y, z)) \, \mid \, x, y, z \in \mathbb{F}_q\big\} \cup \{P_\infty\}\nonumber,
\end{align*}
is an ovoid of \[
Q: X_0X_7 + X_1X_6 + X_2X_5 + X_3X_4 = 0,
\] and $f_1,f_2,f_3\in \mathbb{F}_q[X,Y,Z]$ have degree $2$, then $q$ is even and $\mO_7(f_1,f_2,f_3)$ is the Kantor ovoid as in \eqref{KO}.
\end{theorem}

\section{Degree 3 Case for $f_1, f_2, f_3$}\label{classdeg3}

Throughout this section, we assume that $f_1, f_2, f_3$ are polynomials of degree 3, i.e., 
\begin{eqnarray}\label{Eq:f_i3}
f_1(X,Y,Z) &=& \sum_{i+j+k \leq 3} a_{i,j,k} X^i Y^j Z^k, \nonumber\\
f_2(X,Y,Z) &=& \sum_{i+j+k \leq 3} b_{i,j,k} X^i Y^j Z^k, \\
f_3(X,Y,Z) &=& \sum_{i+j+k \leq 3} c_{i,j,k} X^i Y^j Z^k.\nonumber
\end{eqnarray}

By Theorem \ref{Th:Main}, if $\mO_7(f_1, f_2, f_3)$ is an ovoid of $Q$ and $q > 6.3\cdot 4^{13/3}$ then the hypersurface $\mathcal{S}_{f_1, f_2, f_3}$ does not contain any absolutely irreducible component defined over $\mathbb{F}_q$. Since $\mathcal{S}_{f_1, f_2, f_3}$ is of degree 4, it must split into either four hyperplanes or two quadrics (interchanged by the Frobenius automorphism $\varphi_q$ defined above). In the following sections, we treat these two cases separately.

\subsection{Case of four hyperplanes}\label{4iperpiani}

In this section, we assume that $\mathcal{S}_{f_1, f_2, f_3}$ splits into four hyperplanes. Arguing as in Section \ref{sec:classdeg2}, the four hyperplanes must be defined by $\pi_i = 0$, for $i \in \{1, 2, 3, 4\}$, where
\begin{eqnarray*}
\pi_1 &=& (X_1 - X_4) + \alpha_1 (X_2 - X_5) + \beta_1 (X_3 - X_6), \\
\pi_2 &=& (X_1 - X_4) + \alpha_2 (X_2 - X_5) + \beta_2 (X_3 - X_6), \\
\pi_3 &=& (X_1 - X_4) + \alpha_3 (X_2 - X_5) + \beta_3 (X_3 - X_6), \\
\pi_4 &=& (X_1 - X_4) + \alpha_4 (X_2 - X_5) + \beta_4 (X_3 - X_6).
\end{eqnarray*}

Therefore, all coefficients of the polynomial 
$$G(1, X_1, X_2, X_3, X_4, X_5, X_6) := F(1, X_1, X_2, X_3, X_4, X_5, X_6) - \pi_1 \pi_2 \pi_3 \pi_4$$
must vanish. Note that $\{1, \alpha_i, \beta_i\}$ must be $\mathbb{F}_q$-linearly independent, otherwise $\mO_7(f_1, f_2, f_3)$ would not be an ovoid. Moreover, since $\mathcal{S}_{f_1, f_2, f_3}$ is defined over $\mathbb{F}_q$, one of the following two cases holds:
\begin{itemize}
    \item $\alpha_i = \alpha^{q^{i-1}}$ and $\beta_i = \beta^{q^{i-1}}$ for some $\alpha, \beta \in \mathbb{F}_{q^4}$; or
    \item $\alpha_2 = \alpha_1^q, \ \beta_2 = \beta_1^q, \ \alpha_4 = \alpha_3^q, \ \beta_4 = \beta_3^q$ for some $\alpha_1, \beta_1, \alpha_3, \beta_3 \in \mathbb{F}_{q^2}$.
\end{itemize}

It is easily seen that the latter case cannot occur, since $\{1, \alpha_i, \beta_i\} \subseteq \mathbb{F}_{q^2}$ cannot be $\mathbb{F}_q$-linearly independent. We will now show that the former case also cannot occur.

First, observe that $q$ must be a power of 3. Indeed, the coefficients of $X_1^4$, $X_1^3 X_4$, and $X_1^2 X_2^2$ in $G(1, X_1, X_2, X_3, X_4, X_5, X_6)$ are $-c_{0,0,0} - 1$, $c_{0,0,0} + 4$, and $6$, respectively, and they can vanish simultaneously only if $c_{0,0,0} = -1$ and $q$ is a power of 3. By imposing that all other coefficients of $G(1, X_1, X_2, X_3, X_4, X_5, X_6)$ are zero (see Appendix for the complete list), we obtain a system of equations from which we can express $a_{i,j,k}$, $b_{i,j,k}$, and $c_{i,j,k}$ in terms of $\alpha$ and $\beta$. Specifically, we obtain the following expressions for $f_1$, $f_2$, and $f_3$:
\begin{eqnarray*}
f_1(X, Y, Z) \!&=&\! -\mathrm{Tr}_{q^4/q}(\beta) X^3 - \beta^{q^3 + q^2 + q + 1} Z^3 - \mathrm{Tr}_{q^4/q}(\alpha^{q^3 + q^2 + q} \beta) Y^3 + a_{1,0,0} X + a_{0,1,0} Y, \\
f_2(X, Y, Z) \!&=&\! -\mathrm{Tr}_{q^4/q}(\alpha) X^3 - \alpha^{q^3 + q^2 + q + 1} Y^3 - \mathrm{Tr}_{q^4/q}(\beta^{q^3 + q^2 + q} \alpha) Z^3 + b_{1,0,0} X - a_{0,1,0} Z, \\
f_3(X, Y, Z)\! &=&\! -X^3 - \mathrm{Tr}_{q^4/q}(\alpha^{q^2 + q + 1}) Y^3 - \mathrm{Tr}_{q^4/q}(\beta^{q^2 + q + 1}) Y^3 - b_{1,0,0} Y - a_{1,0,0} Z.
\end{eqnarray*}

These expressions hold under the following conditions on $\alpha$ and $\beta$:
\begin{eqnarray*}
\mathrm{Tr}_{q^4/q}(\alpha^{q+1}) + \mathrm{Tr}_{q^4/q^2}(\alpha^{q^2+1}) &=& 0, \\
\mathrm{Tr}_{q^4/q}(\beta^{q+1}) + \mathrm{Tr}_{q^4/q^2}(\beta^{q^2+1}) &=& 0, \\
\mathrm{Tr}_{q^4/q}(\alpha^{q+1} \beta^{q^3 + q^2}) + \mathrm{Tr}_{q^4/q^2}(\alpha^{q^2+1} \beta^{q^3 + q}) &=& 0, \\
\mathrm{Tr}_{q^4/q}(\alpha \beta^q + \alpha \beta^{q^2} + \alpha \beta^{q^3}) &=& 0, \\
\mathrm{Tr}_{q^4/q}(\beta \alpha^{q^2+q} + \beta \alpha^{q^3+q} + \beta \alpha^{q^3+q^2}) &=& 0, \\
\mathrm{Tr}_{q^4/q}(\alpha \beta^{q^2+q} + \alpha \beta^{q^3+q} + \alpha \beta^{q^3+q^2}) &=& 0.
\end{eqnarray*}

Now, we show that there are no $\alpha, \beta \in \mathbb{F}_{q^4}$ that satisfy these conditions while keeping $\{1, \alpha, \beta\}$ $\mathbb{F}_q$-linearly independent. First, observe that if $\alpha + \alpha^q + \alpha^{q^2} = 0$, then $\alpha \in \mathbb{F}_{q^4} \cap \mathbb{F}_{q^3} = \mathbb{F}_q$, and therefore $\{1, \alpha\}$ are $\mathbb{F}_q$-linearly dependent. Similarly, $\beta + \beta^q + \beta^{q^2} \neq 0$. From the equation
$$
\mathrm{Tr}_{q^4/q}(\alpha^{q+1}) + \mathrm{Tr}_{q^4/q^2}(\alpha^{q^2+1}) = 0 = \mathrm{Tr}_{q^4/q}(\beta^{q+1}) + \mathrm{Tr}_{q^4/q^2}(\beta^{q^2+1}),
$$
we can write
\begin{eqnarray}\label{A4}
\alpha^{q^3} &=& -\frac{\alpha^{q+1} + \alpha^{q^2+1} + \alpha^{q^2+q}}{\alpha + \alpha^q + \alpha^{q^2}}, \nonumber\\
\beta^{q^3} &=& -\frac{\beta^{q+1} + \beta^{q^2+1} + \beta^{q^2+q}}{\beta + \beta^q + \beta^{q^2}}.
\end{eqnarray}
Substituting these into the equation $\mathrm{Tr}_{q^4/q}(\alpha \beta^q + \alpha \beta^{q^2} + \alpha \beta^{q^3}) = 0$, we obtain
$$
\alpha \beta^q + \alpha^q \beta^{q^2} + \alpha^{q^2} \beta - \alpha \beta^{q^2} - \alpha^q \beta - \alpha^{q^2} \beta^q = 0.
$$
Thus, $\beta$ is a root of the linearized polynomial
$$
(\alpha^q - \alpha) X^{q^2} + (\alpha - \alpha^{q^2}) X^q + (\alpha^{q^2} - \alpha^q) X = 0.
$$
It can be shown that the solutions of this polynomial are of the form $\lambda_1 + \lambda_2 \alpha$ for $\lambda_1, \lambda_2 \in \mathbb{F}_q$, implying that $\{1, \alpha, \beta\}$ are not $\mathbb{F}_q$-linearly independent, leading to a contradiction.

Thus, we have proven the following theorem:

\begin{theorem}\label{Th:d=3,iper}
Assume $q>6.3\cdot 4^{13/3}$ and $f_1,f_2,f_3\in \mathbb{F}_q[X,Y,Z]$ of degree $3$. If $\mathcal{S}_{f_1, f_2, f_3}: F(X_0,X_1,X_2,X_3,X_4,X_5,X_6)=0$ as in \eqref{Eq:F}, splits into four hyperplanes, then the set
\begin{align*}
\label{Eq:Param}
\mO_7(f_1,f_2,f_3):=&\big\{(1, x, y, z, f_1(x, y, z), f_2(x, y, z), f_3(x, y, z),\\
&\hspace*{1.5 cm}-z f_1(x, y, z) - y f_2(x, y, z) - x f_3(x, y, z)) \, \mid \, x, y, z \in \mathbb{F}_q\big\} \cup \{P_\infty\}\nonumber,
\end{align*}
is not an ovoid of \[
Q: X_0X_7 + X_1X_6 + X_2X_5 + X_3X_4 = 0.
\] 
\end{theorem}

\subsection{Case of two quadrics}\label{twoquadrics}
In this section, let \( q \) be a power of a prime \( p \) and assume that \( \mathcal{S}_{f_1, f_2, f_3} \) splits into two quadrics switched by the Frobenius \( \varphi_q \). Let \( \xi \in \mathbb{F}_{q^2} \setminus \mathbb{F}_q \) be such that \( \{ 1, \xi \} \) is a basis of \( \mathbb{F}_{q^2} \) over \( \mathbb{F}_q \).

Then, the two factors of \( F(X_0, X_1, X_2, X_3, X_4, X_5, X_6) \) can be written as
\[
R(X_0, X_1, X_2, X_3, X_4, X_5, X_6) + \xi S(X_0, X_1, X_2, X_3, X_4, X_5, X_6),
\]
\[
R(X_0, X_1, X_2, X_3, X_4, X_5, X_6) + \xi^q S(X_0, X_1, X_2, X_3, X_4, X_5, X_6),
\]
where \( R \) and \( S \) are degree-2 polynomials defined over \( \mathbb{F}_q \).  Thus, affine \( \mathbb{F}_q \)-rational points of \( \mathcal{S}_{f_1, f_2, f_3} \) satisfy
\[
\begin{cases}
R(1, X_1, X_2, X_3, X_4, X_5, X_6) = 0, \\
S(1, X_1, X_2, X_3, X_4, X_5, X_6) = 0.
\end{cases}
\]
Let \( U = X_1 - X_4 \), \( V = X_2 - X_5 \), \( W = X_3 - X_6 \). Note that if $S$ is the zero polynomial then $R=0$ possesses solutions beyond $U=V=Z=0$. Thus, we will assume $S\not\equiv 0$. 

As the polynomials \( R \) and \( S \) must vanish on \( U = V = W = 0 \), they can be written as
\[
Q_R(U, V, W) + U L_R(X_4, X_5, X_6) + V M_R(X_4, X_5, X_6) + W N_R(X_4, X_5, X_6),
\]
and
\[
Q_S(U, V, W) + U L_S(X_4, X_5, X_6) + V M_S(X_4, X_5, X_6) + W N_S(X_4, X_5, X_6),
\]
respectively, where the polynomials \( L_R, L_S, M_R, M_S, N_R, N_S \) are linear in \( X_4, X_5, X_6 \) and \( Q_R \), \( Q_S \) are homogeneous of degree 2.

Suppose that \( (L_S, M_S, N_S) \neq \lambda (L_R, M_R, N_R) \) (as polynomials) for all \( \lambda \in \mathbb{F}_q \). Then there exist instances \( (\overline{x}_4, \overline{x}_5, \overline{x}_6) \in \mathbb{F}_q^3 \) such that
{\small \[
(L_S(\overline{x}_4, \overline{x}_5, \overline{x}_6), M_S(\overline{x}_4, \overline{x}_5, \overline{x}_6), N_S(\overline{x}_4, \overline{x}_5, \overline{x}_6)) \neq \lambda (L_R(\overline{x}_4, \overline{x}_5, \overline{x}_6), M_R(\overline{x}_4, \overline{x}_5, \overline{x}_6), N_R(\overline{x}_4, \overline{x}_5, \overline{x}_6)),
\]}
for all \( \lambda \in \mathbb{F}_q \) (so also non-vanishing), and thus the origin \( (0, 0, 0) \) is a simple \( \mathbb{F}_q \)-rational point for the curve
\[
\begin{cases}
R(1, U, V, W, \overline{x}_4, \overline{x}_5, \overline{x}_6) = 0, \\
S(1, U, V, W, \overline{x}_4, \overline{x}_5, \overline{x}_6) = 0.
\end{cases}
\]
This means that this curve contains an \( \mathbb{F}_q \)-rational component passing through \( (0, 0, 0) \), and this proves the existence of solutions \( (\overline{x}_1, \overline{x}_2, \overline{x}_3, \overline{x}_4, \overline{x}_5, \overline{x}_6) \in \mathbb{F}_q^6 \) to
\begin{equation}\label{sistemaRS}
\begin{cases}
R(1, X_1, X_2, X_3, X_4, X_5, X_6) = 0, \\
S(1, X_1, X_2, X_3, X_4, X_5, X_6) = 0,
\end{cases}
\end{equation}
apart from the trivial ones. Therefore, \( \mathcal{O}_7(f_1, f_2, f_3) \) is not an ovoid in this case.

Thus, we can suppose that there exists \( \lambda \in \mathbb{F}_q \) such that \( (L_S, M_S, N_S) = \lambda (L_R, M_R, N_R) \). Actually, up to taking a suitable linear combination of the equations of System \eqref{sistemaRS}, we may assume \( \lambda = 0 \).

In the following, we use the following notations:
\begin{eqnarray*}
Q_R(U, V, W) &=& A_1 U^2 + A_2 V^2 + A_3 W^2 + A_4 UV + A_5 UW + A_6 VW,\\
Q_S(U, V, W) &=& A_1' U^2 + A_2' V^2 + A_3' W^2 + A_4' UV + A_5' UW + A_6' VW,\\
L_R(X_4, X_5, X_6) &=& B_1 X_4 + B_2 X_5 + B_3 X_6 + B_4,\\
M_R(X_4, X_5, X_6) &=& C_1 X_4 + C_2 X_5 + C_3 X_6 + C_4,\\
N_R(X_4, X_5, X_6) &=& D_1 X_4 + D_2 X_5 + D_3 X_6 + D_4,
\end{eqnarray*}
and we will investigate the cases \( p > 3 \) and \( p = 2 \).

\subsubsection{Case \( p > 3 \)}
In this section, we take \( \xi \in \mathbb{F}_{q^2} \setminus \mathbb{F}_q \) such that \( \xi^2 = k \), where \( k \) is not a square in \( \mathbb{F}_q \). So, we have
\begin{eqnarray*}
F(X_0, X_1, X_2, X_3, X_4, X_5, X_6) &=& R(X_0, X_1, X_2, X_3, X_4, X_5, X_6)^2 \\
&&- k S(X_0, X_1, X_2, X_3, X_4, X_5, X_6)^2,\\
R(X_0, X_1, X_2, X_3, X_4, X_5, X_6) &=& Q_R(U, V, W) + U L_R(X_4, X_5, X_6)\\
&&+ V M_R(X_4, X_5, X_6) + W N_R(X_4, X_5, X_6),\\
S(X_0, X_1, X_2, X_3, X_4, X_5, X_6) &=& Q_S(U, V, W).
\end{eqnarray*}

As \( F - (R^2 - k S^2) \) must be the zero polynomial, we obtain a list of conditions on the coefficients \( A_i, A_i', B_i, C_i, D_i, a_{i,j,k}, b_{i,j,k}, c_{i,j,k} \). Let \( \mathcal{B} \) be the set of all these conditions. Then, all the coefficients \( a_{i,j,k}, b_{i,j,k}, c_{i,j,k} \) can be written in terms of \( A_i, A_i', B_i, C_i, D_i \) (see Appendix).

After substituting these conditions into System \( \mathcal{B} \), we obtain
\[
B_1 \left( A_1 - \frac{B_1}{2} \right) = 0 = A_1^2 - \frac{2}{3} A_1 B_1 - {A_1'}^2 k.
\]
So, either \( B_1 = 0 \) or \( A_1 - \frac{B_1}{2} = 0 \). However, if \( B_1 = 0 \), then from the second equation we have \( A_1^2 - {A_1'}^2 k = 0 \), which implies \( A_1 = A_1' = 0 \) since \( k \) is a non-square in \( \mathbb{F}_q \). Therefore, we may assume in any case that \( A_1 = \frac{B_1}{2} \).

The same argument can be applied to the following pairs of equations in \( \mathcal{B} \):
\[
C_2 \left( A_2 - \frac{C_2}{2} \right) = 0 = A_2^2 - \frac{2}{3} A_2 C_2 - {A_2'}^2 k,
\]
and
\[
D_3 \left( A_3 - \frac{D_3}{2} \right) = 0 = A_3^2 - \frac{2}{3} A_3 D_3 - {A_3'}^2 k,
\]
which yield \( A_2 = \frac{C_2}{2} \) and \( A_3 = \frac{D_3}{2} \), respectively. Substituting the expressions for \( A_1 \), \( A_2 \), and \( A_3 \) into \( \mathcal{B} \) yields
\begin{equation}\label{condizioni1}
{A_1'}^2 k + \frac{1}{12} B_1^2 = {A_2'}^2 k + \frac{1}{12} C_2^2 = {A_3'}^2 k + \frac{1}{12} D_3^2 = 0.
\end{equation}

Assume first that \( q \equiv 1 \pmod{3} \). Then, since \( -3 \) is a square in \( \mathbb{F}_q \) and \( k \) is a non-square in \( \mathbb{F}_q \), we have \( A_1' = B_1 = A_2' = C_2 = A_3' = D_3 = 0 \), which also implies \( A_1 = A_2 = A_3 = 0 \).

Then, from \( \mathcal{B} \) we have the following equations:
\[
{A_4'}^2 k = (A_4 - B_2)^2, \quad {A_5'}^2 k = (A_5 - D_1)^2, \quad {A_6'}^2 k = (A_6 - D_2)^2,
\]
and thus \( A_4' = A_5' = A_6' = 0 \), \( A_4 = B_2 \), \( A_5 = D_1 \), and \( A_6 = D_2 \). Finally, substituting all this information in \( \mathcal{B} \) proves that also \( A_4 = A_5 = A_6 = 0 \), a contradiction to \( R \) and \( S \) being polynomials of degree 2.

\begin{theorem}\label{Th:d=3,quadriche1}
Let $q>6.3\cdot 4^{13/3}$, \( q \equiv 1 \pmod{3} \), $p>3$, $f_1,f_2,f_3\in \mathbb{F}_q[X,Y,Z]$  of degree $3$, and $\mathcal{S}_{f_1, f_2, f_3}$, defined by the polynomial $F(1,X_1,X_2,X_3,X_4,X_5,X_6)$ as in \eqref{Eq:F}, split into two quadrics. The set
\begin{align*}
\label{Eq:Param}
\mO_7(f_1,f_2,f_3):=&\big\{(1, x, y, z, f_1(x, y, z), f_2(x, y, z), f_3(x, y, z),\\
&\hspace*{1.5 cm}-z f_1(x, y, z) - y f_2(x, y, z) - x f_3(x, y, z)) \, \mid \, x, y, z \in \mathbb{F}_q\big\} \cup \{P_\infty\}\nonumber,
\end{align*}
is not an ovoid of \[
Q: X_0X_7 + X_1X_6 + X_2X_5 + X_3X_4 = 0.
\] 
\end{theorem}

We can now assume in the following that \( q \) is odd and \( q \equiv 2 \pmod{3} \). Then, the following are equations in \( \mathcal{B} \) 
\begin{eqnarray*}
B_1 \left( A_4 - \frac{B_2}{2} - \frac{C_1}{2} \right) = 
B_2 \left( A_4 - \frac{B_2}{2} - \frac{C_1}{2} \right) =
B_4 \left( A_4 - \frac{B_2}{2} - \frac{C_1}{2} \right) &=& 0, \\
C_1 \left( A_4 - \frac{B_2}{2} - \frac{C_1}{2} \right) =
C_2 \left( A_4 - \frac{B_2}{2} - \frac{C_1}{2} \right) =
C_4 \left( A_4 - \frac{B_2}{2} - \frac{C_1}{2} \right) &=& 0,
\end{eqnarray*}
and
\[
A_4^2 - 2 A_4 B_2 - 2 A_1' A_2' k - {A_4'}^2 k - \frac{1}{2} B_1 C_2 + B_2^2 = 0.
\]
Assume first that \( A_4 \neq \frac{B_2}{2} + \frac{C_1}{2} \). Then, \( B_1 = B_2 = B_4 = C_1 = C_2 = C_4 = 0 \), which implies \( A_1 = A_2 = A_1' = A_2' = 0 \). Therefore, the equation
\[
A_4^2 - 2 A_4 B_2 - 2 A_1' A_2' k - {A_4'}^2 k - \frac{1}{2} B_1 C_2 + B_2^2 = 0
\]
becomes
\[
A_4^2 - {A_4'}^2 k = 0,
\]
which implies \( A_4 = A_4' = 0 \), contradicting \( A_4 \neq \frac{B_2}{2} + \frac{C_1}{2} = 0 \).

Thus, we may assume
\[
A_4 = \frac{1}{2} \left( B_2 + C_1 \right).
\]

The same argument can be applied to the equations
\begin{eqnarray*}
B_1 \left( A_5 - \frac{B_3}{2} - \frac{D_1}{2} \right) =
B_3 \left( A_5 - \frac{B_3}{2} - \frac{D_1}{2} \right) =
B_4 \left( A_5 - \frac{B_3}{2} - \frac{D_1}{2} \right) &=& 0, \\
D_1 \left( A_5 - \frac{B_3}{2} - \frac{D_1}{2} \right) =
D_3 \left( A_5 - \frac{B_3}{2} - \frac{D_1}{2} \right)=
D_4 \left( A_5 - \frac{B_3}{2} - \frac{D_1}{2} \right) &=& 0,
\end{eqnarray*}
and
\[
A_5^2 - 2 A_5 B_3 - 2 A_1' A_3' k - {A_5'}^2 k - \frac{1}{2} B_1 D_3 + B_3^2 = 0,
\]
and the equations
\begin{eqnarray*}
D_2 \left( A_6 - \frac{D_2}{2} - \frac{C_3}{2} \right) =
D_3 \left( A_6 - \frac{D_2}{2} - \frac{C_3}{2} \right) =
D_4 \left( A_6 - \frac{D_2}{2} - \frac{C_3}{2} \right) &=& 0, \\
C_2 \left( A_6 - \frac{D_2}{2} - \frac{C_3}{2} \right) =
C_3 \left( A_6 - \frac{D_2}{2} - \frac{C_3}{2} \right) =
C_4 \left( A_6 - \frac{D_2}{2} - \frac{C_3}{2} \right) &=& 0,
\end{eqnarray*}
which yield
\[
A_5 = \frac{1}{2} \left( B_3 + D_1 \right), \quad A_6 = \frac{1}{2} \left( D_2 + C_3 \right).
\]

We now distinguish two cases according to whether \( B_1 \) is zero or not.

\begin{itemize}
    \item[(1)] \textbf{CASE \( B_1 \neq 0 \)}. Observe first that by \eqref{condizioni1} and \( B_1 \neq 0 \), we have \( A_1' \neq 0 \), and
  \begin{equation}\label{kcaso1}
    k = -\frac{B_1^2}{12{A_1'}^2}.
  \end{equation}
Now, consider the following equations in \( \mathcal{B} \):
\[
A_1' A_4' k + \frac{1}{12} B_1 B_2 + \frac{1}{12} B_1 C_1 = 0 = A_1' A_5' k + \frac{1}{12} B_1 B_3 + \frac{1}{12} B_1 D_1.
\]
By replacing \eqref{kcaso1} in the above conditions, we obtain
\[
B_1 A_1' \left( A_1' (B_2 + C_1) - B_1 A_4' \right) = 0 = B_1 A_1' \left( A_1' (B_3 + D_1) - B_1 A_5' \right),
\]
which allows us to obtain
\[
A_4' = \frac{A_1'}{B_1} (B_2 + C_1), \quad A_5' = \frac{A_1'}{B_1} (B_3 + D_1).
\]
Similarly, one can obtain
\[
A_6' = \frac{A_1'}{B_1} (C_3 + D_2).
\]
By substituting this information in \( \mathcal{B} \), we obtain the following conditions:
\[
A_2' = \frac{A_1'}{B_1} C_2, \quad A_3' = \frac{A_1'}{B_1} D_3.
\]
These equations allow us to write all the conditions in \( \mathcal{B} \) in terms of \( B_i, C_i, D_i \). After doing so, we obtain the following equations:
\begin{eqnarray*}
C_2 D_3 - C_3^2 + C_3 D_2 - D_2^2 &=& 0,\\
B_2 D_3 - 2 B_3 C_3 + B_3 D_2 + C_1 D_3 + C_3 D_1 - 2 D_1 D_2 &=& 0,\\
B_1 C_2 - B_2^2 + B_2 C_1 - C_1^2 &=& 0,\\
B_2 C_3 - 2 B_2 D_2 + B_3 C_2 - 2 C_1 C_3 + C_1 D_2 + C_2 D_1 &=& 0,\\
B_1 D_3 - B_3^2 + B_3 D_1 - D_1^2 &=& 0,\\
B_1 C_3 + B_1 D_2 - 2 B_2 B_3 + B_2 D_1 + B_3 C_1 - 2 C_1 D_1 &=& 0.
\end{eqnarray*}
Since \( B_1 \neq 0 \), from these equations we have
\begin{eqnarray*}
C_2 &=& \frac{B_2^2 - B_2 C_1 + C_1^2}{B_1},\\
D_3 &=& \frac{B_3^2 - B_3 D_1 + D_1^2}{B_1},\\
D_2 &=& \frac{2 B_2 B_3 - B_1 C_3 - B_2 D_1 - B_3 C_1 + 2 C_1 D_1}{B_1}.
\end{eqnarray*}
After these substitutions, we obtain the following conditions:
\begin{eqnarray*}
A_1'^2 B_1^5 (B_3 - D_1) (B_1 C_3 - B_2 B_3 + B_2 D_1 - C_1 D_1) &=& 0,\\
A_1'^2 B_1^6 (B_1 C_3 - B_2 B_3 + B_3 C_1 - C_1 D_1) (B_1 C_3 - B_2 B_3 + B_2 D_1 - C_1 D_1) &=& 0,\\
A_1'^2 B_1^5 (B_2 - C_1) (B_1 C_3 - B_2 B_3 + B_2 D_1 - C_1 D_1) &=& 0.
\end{eqnarray*}
\end{itemize}

Then, observe that
\[
B_1C_3 - B_2B_3 + B_2D_1 - C_1D_1 = 0.
\]
Indeed, if this is not true, then \(B_3 = D_1\), \(B_2 = C_1\), and
\[
0 = B_1C_3 - B_2B_3 + B_3C_1 - C_1D_1 = B_1C_3 - B_3B_2 = B_1C_3 - B_2B_3 + B_2D_1 - C_1D_1 \neq 0,
\]
which is a contradiction. If \(B_1C_3 - B_2B_3 + B_2D_1 - C_1D_1 = 0\), all the conditions are satisfied and we have no contradiction.

\item[(2)] \textbf{CASE \(B_1 = 0\).} In this case, \({A_1^\prime} = 0\). First, observe that it is not possible that \({A_2^\prime} = 0 = {A_3^\prime}\), otherwise, we could argue as in the case \(q \equiv 1 \pmod{3}\) and obtain a contradiction. So, assume first \({A_2^\prime} \neq 0\), and hence \(C_2 \neq 0\). Then we can use
\[
k = -\frac{C_2^2}{12{A_2^\prime}^2},
\]
and arguing as in \textbf{CASE} (1), we obtain
\[
{A_3^\prime} = \frac{A_2^\prime}{C_2}D_3, \quad {A_4^\prime} = \frac{A_2^\prime}{C_2}(B_2 + C_1), \quad {A_5^\prime} = \frac{A_2^\prime}{C_2}(B_3 + D_1), \quad {A_6^\prime} = \frac{A_2^\prime}{C_2}(C_3 + D_2),
\]
and so the conditions in \(\mathcal{B}\) read
\begin{eqnarray*}
C_2D_3 - C_3^2 + C_3D_2 - D_2^2 &=& 0,\\
B_2D_3 - 2B_3C_3 + B_3D_2 + C_1D_3 + C_3D_1 - 2D_1D_2 &=& 0,\\
B_2^2 - B_2C_1 + C_1^2 &=& 0,\\
B_2C_3 - 2B_2D_2 + B_3C_2 - 2C_1C_3 + C_1D_2 + C_2D_1 &=& 0,\\
B_3^2 - B_3D_1 + D_1^2 &=& 0,\\
-2B_2B_3 + B_2D_1 + B_3C_1 - 2C_1D_1 &=& 0.
\end{eqnarray*}
Observe that since \(q \equiv 2 \pmod{3}\), \(-3\) is a non-square of \(\mathbb{F}_q\), and hence \(B_2^2 - B_2C_1 + C_1^2 = 0\) and \(B_3^2 - B_3D_1 + D_1^2 = 0\) imply \(B_2 = C_1 = B_3 = D_1 = 0\).

So, in this case, all conditions are satisfied whenever
\[
C_2D_3 - C_3^2 + C_3D_2 - D_2^2 = 0.
\]

In this case, the System \eqref{sistemaRS} reads
\[
\begin{cases}
\frac{C_2}{2}V^2 + \frac{D_3}{2}W^2 + \frac{D_2 + C_3}{2}VW + B_4U\\
\hspace*{1 cm}+ V(C_2X_5 + C_3X_6 + C_4) + W(D_2X_5 + D_3X_6 + D_4) = 0, \\
A_2^\prime V^2 + \frac{A_2^\prime}{C_2}D_3 W^2 + \frac{A_2^\prime}{C_2}(D_2 + C_3)VW = 0,
\end{cases}
\]
which simplifies to
\[
\begin{cases}
C_2V^2 + D_3W^2 + (D_2 + C_3)VW + 2B_4U \\
\hspace*{1.5 cm}+ 2V(C_2X_5 + C_3X_6 + C_4)+ 2W(D_2X_5 + D_3X_6 + D_4) = 0, \\
C_2V^2 + D_3W^2 + (D_2 + C_3)VW = 0,
\end{cases}
\]
and hence we obtain the system
\[
\begin{cases}
2B_4U + 2V(C_2X_5 + C_3X_6 + C_4) + 2W(D_2X_5 + D_3X_6 + D_4) = 0, \\
C_2V^2 + D_3W^2 + (D_2 + C_3)VW = 0.
\end{cases}
\]
Now, by the condition on \(C_2, C_3, D_2, D_3\), the second equation factorizes as
\[
\begin{cases}
2B_4U + 2V(C_2X_5 + C_3X_6) + 2W(D_2X_5 + D_3X_6 + D_4) = 0, \\
\left(V + \frac{D_2 + C_3}{2C_2}W + \eta W(C_3 - D_2)\right)\left(V + \frac{D_2 + C_3}{2C_2}W - \eta W(C_3 - D_2)\right) = 0,
\end{cases}
\]
where \(\eta \in \mathbb{F}_{q^2} \setminus \mathbb{F}_q\) is such that \(\eta^2 = -3\).

Clearly, if \(B_4 = 0\), then \((u, 0, 0)\) is a solution for every \(u \in \mathbb{F}_q\), and hence \(\mathcal{O}_7(f_1, f_2, f_3)\) is not an ovoid of \(Q\). Observe that the same holds if \(C_3 = D_2\), in which case the above system reads
\[
\begin{cases}
B_4 U = 0, \\
\left(V + \frac{C_3}{C_2}W\right)^2 = 0,
\end{cases}
\]
and hence it always admits a non-trivial solution \((U, V, W) \in \mathbb{F}_q^3 \setminus \{(0, 0, 0)\}\) and $\mathcal{O}_7(f_1,f_2,f_3)$ is not an ovoid of $Q$. 

In the remaining cases, i.e., if \(B_4 \neq 0\) and \(C_3 \neq D_2\), we have no contradictions.

Finally, assume \({A_2^\prime} = 0 = C_2\), and hence \({A_3^\prime} \neq 0\) and \(D_3 \neq 0\). Then, by repeating the above computations using
\[
k = -\frac{D_3^2}{12{A_3^\prime}^2},
\]
we have
\[
{A_4^\prime} = \frac{A_3^\prime}{D_3}(B_2 + C_1), \quad {A_5^\prime} = \frac{A_3^\prime}{D_3}(B_3 + D_1), \quad {A_6^\prime} = \frac{A_3^\prime}{D_3}(C_3 + D_2),
\]
and the conditions
\begin{eqnarray*}
C_3^2 - C_3D_2 + D_2^2 &=& 0,\\
B_2D_3 - 2B_3C_3 + B_3D_2 + C_1D_3 + C_3D_1 - 2D_1D_2 &=& 0,\\
B_2^2 - B_2C_1 + C_1^2 &=& 0,\\
B_2C_3 - 2B_2D_2 - 2C_1C_3 + C_1D_2 &=& 0,\\
B_3^2 - B_3D_1 + D_1^2 &=& 0,\\
-2B_2B_3 + B_2D_1 + B_3C_1 - 2C_1D_1 &=& 0.
\end{eqnarray*}
As before, this yields \(C_3 = D_2 = B_2 = C_1 = B_3 = D_1 = 0\), and System \eqref{sistemaRS} reads
\[
\begin{cases}
\frac{D_3}{2}W^2 + B_4U + C_4V + W(D_3 X_6 + D_4) = 0, \\
A_3^\prime W^2 = 0.
\end{cases}
\]
Clearly, this system always has solutions \((U, V, W) \neq (0, 0, 0)\), and hence \(\mathcal{O}_7(f_1, f_2, f_3)\) is not an ovoid of \(Q\) in this case.

Observe that in the open cases (where no contradiction was obtained), we always have \(Q_R = \alpha Q_S\), and hence System \ref{sistemaRS} reads
\[
\begin{cases}
\alpha Q_S(U, V, W) + U L_R(X_4, X_5, X_6) + V M_R(X_4, X_5, X_6) + W N_R(X_4, X_5, X_6) = 0, \\
Q_S(U, V, W) = 0,
\end{cases}
\]
which simplifies to
\[
\begin{cases}
U L_R(X_4, X_5, X_6) + V M_R(X_4, X_5, X_6) + W N_R(X_4, X_5, X_6) = 0, \\
Q_S(U, V, W) = 0.
\end{cases}
\]
Now, \(Q_S(U, V, W) = 0\) can be seen as a (possibly reducible) conic in \(\mathrm{PG}(2, q)\), and hence it always contains a point \((\bar{u} : \bar{v} : \bar{w}) \in \mathrm{PG}(2, q)\). Moreover, we can always choose \((\bar{x}_4, \bar{x}_5, \bar{x}_6) \in \mathbb{F}_q^3\) such that
\[
\bar{u} L_R(\bar{x}_4, \bar{x}_5, \bar{x}_6) + \bar{v} M_R(\bar{x}_4, \bar{x}_5, \bar{x}_6) + \bar{w} N_R(\bar{x}_4, \bar{x}_5, \bar{x}_6) = 0,
\]
and hence the 6-tuple \((\bar{u} + \bar{x}_4, \bar{v} + \bar{x}_5, \bar{w} + \bar{x}_6, \bar{x}_4, \bar{x}_5, \bar{x}_6)\) is a non-trivial solution of the above system unless, up to projectivities, we are in the following situation
\[
\begin{cases}
U + V M_R(X_4, X_5, X_6) + W N_R(X_4, X_5, X_6) = 0, \\
(V + \xi W)(V - \xi W) = 0.
\end{cases}
\]
In this case, the two factors of \(F\) would be
\[
(1 + \xi)(V^2 - kW^2) + U + VM_R + WN_R
\]
and
\[
(1 - \xi)(V^2 - kW^2).
\]
If we now impose this factorization, it is readily seen that this implies \(k = -\frac{1}{3}\), \(C_1 = D_1 = 0\), \(C_2 = 2\), \(C_3 = -D_2 = \epsilon \frac{2}{3}\), and \(D_3 = \frac{2}{3}\), where \(\epsilon = \pm 1\).

Thus, we have the following. 
\begin{theorem}\label{Th:d=3,quadriche2}
Let $q>6.3\cdot 4^{13/3}$, $q\equiv 2 \pmod{3}$, $p>3$, and $f_1,f_2,f_3\in \mathbb{F}_q[X,Y,Z]$ of degree $3$. Assume that $\mathcal{S}_{f_1, f_2, f_3}: F(1,X_1,X_2,X_3,X_4,X_5,X_6)=0$ as in \eqref{Eq:F} splis into two quadrics. 
If \begin{align*}
\label{Eq:Param}
\mO_7(f_1,f_2,f_3):=&\big\{(1, x, y, z, f_1(x, y, z), f_2(x, y, z), f_3(x, y, z),\\
&\hspace*{1.5 cm}-z f_1(x, y, z) - y f_2(x, y, z) - x f_3(x, y, z)) \, \mid \, x, y, z \in \mathbb{F}_q\big\} \cup \{P_\infty\}\nonumber,
\end{align*}
is  an ovoid of \[
Q: X_0X_7 + X_1X_6 + X_2X_5 + X_3X_4 = 0,
\]  
then 
\begin{eqnarray}\label{famiglia1}
f_1(X,Y,Z)&=& -\frac{4}{27}Z^3+\frac{4\epsilon}{3} Y^3-\frac{4}{9}Y^2Z+\frac{4\epsilon}{9}YZ^2+\left(\frac{4\epsilon}{3}C_4-2D_4\right)Y^2-\frac{2}{3}D_4Z^2+\nonumber\\&&+\frac{4\epsilon}{3}XY+\frac{4\epsilon}{3}D_4YZ+a_{0,1,0}Y-D_4^2Z,\nonumber\\
f_2(X,Y,Z)&=& -\frac{4}{3} Y^3-\frac{4\epsilon}{9}Z^3-\frac{4}{9}YZ^2+\frac{4\epsilon}{3}Y^2Z-2C_4Y^2-\left(\frac{4\epsilon}{3}D_4-\frac{2}{3}C_4\right)Z^2\\&&-\frac{4\epsilon}{3}XZ+\frac{4\epsilon}{3}C_4YZ+b_{1,0,0}X-C_4^2Y-(2C_4D_4+a_{0,1,0})Z,\nonumber\\
f_3(X,Y,Z)&=& -2Y^2-\frac{2}{3}Z^2-X-(2C_4+b_{1,0,0})Y-(2D_4+a_{1,0,0})Z.\nonumber
\end{eqnarray} 
\end{theorem}

Section \ref{Kantorq=2} shows that the Kantor ovoid for $q\equiv 2 \pmod{3}$ has degree $3$ and it falls within the family $\mO_7(f_1,f_2,f_3)$ as in \eqref{famiglia1}. We leave as an open problem to determine whether the Kantor ovoid is the unique (up to equivalences) ovoid of degree $3$ of low-degree.

\begin{open}
Determine whether for any $f_1, f_2, f_3$ as in \eqref{famiglia1}, the ovoid $\mO_7(f_1,f_2,f_3)$ is equivalent to the Kantor ovoid as in Section \ref{Kantorq=2}.
\end{open}

\subsubsection{Case $p=2$}
In this section, we will continue with the same notations and assumptions as in the previous sections, and we will take $\xi \in \mathbb{F}_{q^2} \setminus \mathbb{F}_q$ such that $\xi^q = \xi + 1$. 

Arguing as before, imposing that
$$
F - (R + \xi S)(R + \xi^q S)
$$
is the zero polynomial, we obtain a list of conditions $\mathcal{B}$ (see Appendix). Among these conditions, we find
\begin{eqnarray}\label{sistemap=2}
    A_3 A_6^\prime + A_6 A_3^\prime = A_1 A_4^\prime + A_4 A_1^\prime = A_2 A_6^\prime + A_6 A_2^\prime &=& 0, \nonumber\\
    A_1 A_5^\prime + A_5 A_1^\prime = A_3 A_5^\prime + A_5 A_3^\prime = A_2 A_4^\prime + A_4 A_2^\prime &=& 0, \nonumber\\
    A_3 A_4^\prime + A_4 A_3^\prime + A_5 A_6^\prime + A_6 A_5^\prime &=& 0, \nonumber\\
    A_1 A_6^\prime + A_4 A_5^\prime + A_5 A_4^\prime + A_6 A_1^\prime &=& 0, \nonumber\\
    A_2 A_5^\prime + A_4 A_6^\prime + A_5 A_2^\prime + A_6 A_4^\prime &=& 0.
\end{eqnarray}

Using these conditions, our first aim is to prove that 
$(A_1, \ldots, A_6) = \alpha (A_1^\prime, \ldots, A_6^\prime)$ for some $\alpha \in \mathbb{F}_q$. 
Clearly, if $(A_1^\prime, \ldots, A_6^\prime) = (0, \ldots, 0)$, then $S = Q_S$ is the zero polynomial, which contradicts $S$ being a polynomial of degree 2. 
So, without loss of generality, assume $A_1^\prime \neq 0$. Then, from the above equations, we have
$$
A_4 = \frac{A_1}{A_1^\prime} A_4^\prime \coloneqq \alpha A_4^\prime, \quad
A_5 = \frac{A_1}{A_1^\prime} A_5^\prime \coloneqq \alpha A_5^\prime, \quad 
A_6 = \frac{A_1}{A_1^\prime} A_6^\prime \coloneqq \alpha A_6^\prime.
$$
The other conditions in \eqref{sistemap=2} now yield
\begin{eqnarray*}
A_3 A_4^\prime &=& \alpha A_4^\prime A_3^\prime, \quad
A_3 A_5^\prime = \alpha A_5^\prime A_3^\prime, \quad
A_3 A_6^\prime = \alpha A_6^\prime A_3^\prime, \\
A_2 A_4^\prime &=& \alpha A_4^\prime A_2^\prime, \quad
A_2 A_5^\prime = \alpha A_5^\prime A_2^\prime, \quad 
A_2 A_6^\prime = \alpha A_6^\prime A_2^\prime.
\end{eqnarray*}
So, if $(A_4^\prime, A_5^\prime, A_6^\prime) \neq (0,0,0)$, then
$A_2 = \alpha A_2^\prime$ and $A_3 = \alpha A_3^\prime$, and hence 
$(A_1, \ldots, A_6) = \alpha (A_1^\prime, \ldots, A_6^\prime)$ follows.

Therefore, assume $(A_4^\prime, A_5^\prime, A_6^\prime) = (0,0,0)$. Then, $(A_4, A_5, A_6) = (0,0,0)$ holds. 
By substituting this information into $\mathcal{B}$, we obtain $C_3 = D_2$, $B_2 = C_1$, and $B_3 = D_1$, which then imply
$$
A_1 A_3^\prime + A_3 A_1^\prime = 0, \quad \text{and} \quad A_1 A_2^\prime + A_2 A_1^\prime = 0.
$$
So, also in this case $(A_1, \ldots, A_6) = \alpha (A_1^\prime, \ldots, A_6^\prime)$ holds with $\alpha = A_1 / A_1^{\prime}$.

Therefore, System \eqref{sistemaRS} now reads
\begin{equation*}
\begin{cases}
\alpha Q_S(U,V,W) + U L_R(X_4, X_5, X_6) + V M_R(X_4, X_5, X_6) + W N_R(X_4, X_5, X_6) = 0, \\
Q_S(U,V,W) = 0,
\end{cases}
\end{equation*}
that is,
\begin{equation*}
\begin{cases}
U L_R(X_4, X_5, X_6) + V M_R(X_4, X_5, X_6) + W N_R(X_4, X_5, X_6) = 0, \\
Q_S(U,V,W) = 0.
\end{cases}
\end{equation*}

To conclude this section, we argue as in the last part of the case $p > 3$. Namely, $Q_S(U, V, W) = 0$ is a (possibly reducible) conic in $\mathrm{PG}(2, q)$, and hence it always contains a point $(\bar{u} : \bar{v} : \bar{w}) \in \mathrm{PG}(2, q)$. Therefore, $\mO_7(f_1, f_2, f_3)$ is not an ovoid of $Q$ unless, up to projectivities, we are in the following situation
\begin{equation*}
\begin{cases}
U + V M_R(X_4, X_5, X_6) + W N_R(X_4, X_5, X_6) = 0, \\
(V + \xi W)(V + \xi^q W) = 0.
\end{cases}
\end{equation*}

Now, by imposing that 
$$
(1 + \xi)(V + \xi W)(V + \xi^q W) + U + V M_R + W N_R
$$
and
$$
(1 + \xi^q)(V + \xi W)(V + \xi^q W)
$$
are factors of $F$, we obtain $C_1 = D_1 = C_3 = 0$, $C_2 = D_2 = D_3 = 1$.

Thus, we have the following.

\begin{theorem}\label{Th:d=3,quadriche3}
Let $q > 6.3 \cdot 4^{13/3}$, $q \equiv 2 \pmod{3}$, $p = 2$, and $f_1,f_2,f_3\in \mathbb{F}_q[X,Y,Z]$  of degree $3$. Assume that $\mathcal{S}_{f_1, f_2, f_3}: F(1,X_1,X_2,X_3,X_4,X_5,X_6)=0$ as in \eqref{Eq:F} splits into two quadrics.  If \begin{align*}
\label{Eq:Param}
\mO_7(f_1,f_2,f_3):=&\big\{(1, x, y, z, f_1(x, y, z), f_2(x, y, z), f_3(x, y, z),\\
&\hspace*{1.5 cm}-z f_1(x, y, z) - y f_2(x, y, z) - x f_3(x, y, z)) \, \mid \, x, y, z \in \mathbb{F}_q\big\} \cup \{P_\infty\}\nonumber,
\end{align*}
is  an ovoid of \[
Q: X_0X_7 + X_1X_6 + X_2X_5 + X_3X_4 = 0,
\]  
then 
\begin{eqnarray*}\label{famiglia2}
f_1(X,Y,Z) &=& Z^3 + Y^2 Z + Y Z^2 + (C_4 + D_4) Y^2 + D_4 Z^2 + X Y \\
&& + c_{0,0,1} X + b_{0,0,1} Y + D_4^2 Z, \nonumber \\
f_2(X,Y,Z) &=& Y^3 + Z^3 + C_4 Y^2 + (C_4 + D_4) Z^2 + X Z + c_{0,1,0} X + C_4^2 Y + b_{0,0,1} Z, \nonumber \\
f_3(X,Y,Z) &=& Y^2 + Z^2 + Y Z + X + c_{0,1,0} Y + c_{0,0,1} Z. \nonumber
\end{eqnarray*}
\end{theorem}

\begin{open}
Determine whether for $f_1, f_2, f_3$ as in \eqref{famiglia2}, the ovoid $\mO_7(f_1, f_2, f_3)$ is equivalent to the Kantor ovoid as in Section \ref{Kantorq=2}.
\end{open}

\subsubsection{Case $p=3$}
When \( p = 3 \), we can proceed similarly to the case \( p > 3 \) and derive some partial results regarding the functions \( f_1, f_2, f_3 \). However, obtaining a complete classification for \( \mathcal{O}_7(f_1, f_2, f_3) \) in this scenario proved to be challenging, and we leave this as an open problem.

\begin{open}
Obtain a full classification for ovoids $\mO_7(f_1,f_2,f_3)$ of degree $3$ when $p=3$.
\end{open}

\section*{Acknowledgements}
The authors thank the Italian National Group for Algebraic and Geometric Structures and their Applications (GNSAGA—INdAM)
which supported the research. 

\section*{Declarations}
{\bf Conflicts of interest.} The authors have no conflicts of interest to declare that are relevant to the content of this
article.

\appendix
\section{Appendix}

In this Appendix, we give the list of conditions on the coefficients $A_i,B_i,C_i,D_i,E_i,F_i$, see Section \ref{sec:classdeg2}, and $A_i,A_i',B_i,C_i,D_i, a_{i,j,k}, b_{i,j,k}, c_{i,j,k}$, see Subsection \ref{twoquadrics}.

\subsection{Conditions in Section $4$}

\begin{eqnarray*}
    E_2 + \mathrm{N}_{q^3/q}(\alpha)=0,\\
     C_2 + E_1 +\mathrm{Tr}_{q^3/q}(\alpha^{1+q}\beta^{q^2}) =0,\\
        E_3 +\mathrm{Tr}_{q^3/q}(\alpha^{1+q}) =0,\\
        E_2 + 3\mathrm{N}_{q^3/q}(\alpha)=0, \\
        E_1 + \mathrm{Tr}_{q^3/q}(\alpha^{1+q}\beta^{q^2})=0,\\
        C_1 +F_2 +\mathrm{Tr}_{q^3/q}(\alpha \beta^{q+q^2}) =0,\\
        C_3 +\mathrm{Tr}_{q^3/q}(\alpha \beta^q + \alpha^q\beta) =0,\\
        C_2 + 2\mathrm{Tr}_{q^3/q}(\alpha^{1+q}\beta^{q^2})=0,\\
       C_1 + 2\mathrm{Tr}_{q^3/q}(\alpha \beta^{q+q^2})=0,\\
       D_2 +\mathrm{Tr}_{q^3/q}(\alpha)=0,\\
        A_2 - 2\mathrm{Tr}_{q^3/q}(\alpha^{1+q})=0,\\
        B_2 - \mathrm{Tr}_{q^3/q}(\alpha \beta^q+\alpha^q \beta)=0,\\
        E_2 - 3\mathrm{N}_{q^3/q}(\alpha)=0,\\
        C_2 -
        2\mathrm{Tr}_{q^3/q}(\alpha^{1+q}\beta^{q^2})=0,\\
       F_2 - \mathrm{Tr}_{q^3/q}(\alpha \beta^{q+q^2})=0,\\
        F_1 + \mathrm{N}_{q^3/q}(\beta)=0,\\
        F_3 + \mathrm{Tr}_{q^3/q}(\beta^{1+q})=0,\\
        F_2 +\mathrm{Tr}_{q^3/q}(\alpha \beta^{q+q^2}) =0,\\
        F_1 + 3\mathrm{N}_{q^3/q}(\beta)=0,\\
        D_1 -
        \mathrm{Tr}_{q^3/q}(\beta)=0,\\
        A_1 -\mathrm{Tr}_{q^3/q}(\alpha \beta^q+\alpha^q \beta)=0,\\
        B_1 - 2\mathrm{Tr}_{q^3/q}(\beta^{1+q})=0,\\
        E_1 -\mathrm{Tr}_{q^3/q}(\alpha^{1+q}\beta^{q^2}) =0,\\
       C_1 - 2\mathrm{Tr}_{q^3/q}(\alpha \beta^{q+q^2})=0,\\
        F_1 -
        3\mathrm{N}_{q^3/q}(\beta)=0,\\
         D_3 - 1=0,\\
         A_3 + D_2 -
        \mathrm{Tr}_{q^3/q}(\alpha)=0,\\
         B_3 + D_1 -
        \mathrm{Tr}_{q^3/q}(\beta)=0,\\
         A_2 + E_3 -\mathrm{Tr}_{q^3/q}(
        \alpha^{1+q})=0,\\
         A_1 +
        B_2 + C_3 - \mathrm{Tr}_{q^3/q}(\alpha \beta^q+\alpha^q \beta) =0,\\
\end{eqnarray*}

\begin{eqnarray*}
B_1 + F_3 -\mathrm{Tr}_{q^3/q}(\beta^{1+q})=0,\\
       E_2 - \mathrm{N}_{q^3/q}(\alpha)=0,\\
         C_2 + E_1 -\mathrm{Tr}_{q^3/q}( \alpha^{1+q}\beta^{q^2})=0,\\
    C_1 + F_2 -\mathrm{Tr}_{q^3/q}(
        \alpha \beta^{q+q^2})=0,\\
        F_1 -\mathrm{N}_{q^3/q}(\beta)=0,\\
         A_2 + E_3 + \mathrm{Tr}_{q^3/q}(\alpha^{1+q})=0,\\
     A_1 +
        B_2 + C_3 + \mathrm{Tr}_{q^3/q}(\alpha \beta^q+\alpha^q\beta)=0, \\
        A_3 + 2\mathrm{Tr}_{q^3/q}(\alpha)=0, \\
        A_2 + 2\mathrm{Tr}_{q^3/q}(\alpha^{1+q}) =0, \\
        A_1 + \mathrm{Tr}_{q^3/q}(\alpha\beta^q+\alpha^q\beta)=0, \\
         B_1 + F_3 +\mathrm{Tr}_{q^3/q}(\beta^{1+q}) =0, \\
        B_3 + 2\mathrm{Tr}_{q^3/q}(\beta)=0, \\
        B_2 + \mathrm{Tr}_{q^3/q}(\alpha \beta^q+\alpha^q\beta)=0, \\
        B_1 + 2\mathrm{Tr}_{q^3/q}(\beta^{1+q})=0, \\
        D_3 - 3=0, \\
        A_3 - 2\mathrm{Tr}_{q^3/q}(\alpha)=0, \\
        B_3 - 2\mathrm{Tr}_{q^3/q}(\beta)=0, \\
        E_3 - \mathrm{Tr}_{q^3/q}(\alpha^{1+q})=0, \\
        C_3 -\mathrm{Tr}_{q^3/q}(\alpha \beta^q+\alpha^q\beta)=0, \\
        F_3 - \mathrm{Tr}(\beta^{1+q})=0, \\
         A_3 + D_2 + \mathrm{Tr}_{q^3/q}(\alpha)=0, \\
         B_3 +D_1+ \mathrm{Tr}_{q^3/q}(\beta)=0,\\
        D_3 + 3=0, \\
        D_2 + \mathrm{Tr}_{q^3/q}(\alpha)=0, \\
        D_1 + \mathrm{Tr}_{q^3/q}(\beta)=0, \\
        D_3+1=0.\\
\end{eqnarray*}

\subsection{Conditions in Subsection $5.2$}

\textbf{CASE $p>2$}

\begin{eqnarray*}
    B_4^2+c_{1,0,0}=0, \\ 
    2A_1B_4+c_{2,0,0}=0, \\ 
    A_1^2-k{A_1'}^2+c_{3,0,0}=0, \\ 
    C_4^2+b_{0,1,0}=0, \\ 
    2A_2C_4+b_{0,2,0}=0, \\ 
    A_2^2-k{A_2'}^2+b_{0,3,0}=0, \\ 
    D_4^2+a_{0,0,1}=0, \\ 
    2A_3D_4+a_{0,0,2}=0, \\ 
    A_3^2-k{A_3'}^2+a_{0,0,3}=0, \\ 
    -2A_1B_4 +2B_1B_4+c_{2,0,0}=0, \\ 
    A_1^2-2A_1B_1-{A_1'}^2+B_1^2+c_{3,0,0}=0,\\ 
    -2A_2C_4+2C_2C_4+b_{0,2,0}=0, \\ 
    A_2^2-2A_2C_2-k{A_2'}^2+C_2^2+b_{0,3,0}=0, \\ 
    -2A_3D_4+2D_3D_4+a_{0,0,2}=0, \\ 
    A_3^2-2A_3D_3-k{A_3'}^2+D_3^2+a_{0,0,3}=0, \\ 
     2B_4C_4+b_{1,0,0}+c_{0,1,0}=0, \\ 
    2A_2B_4+2A_4C_4+b_{1,1,0}+c_{0,2,0}=0, \\ 
    2A_2A_4-2kA_2'A_4'+b_{1,2,0}+c_{0,3,0}=0, \\ 
2A_1C_4+2A_4B_4+b_{2,0,0}+c_{1,1,0}=0, \\ 
    2A_1A_2+A_4^2-2kA_1'A_2'-k{A_4'}^2+b_{2,1,0}+c_{1,2,0}=0, \\ 
    2A_1A_4-2kA_1'A_4'+b_{3,0,0}+c_{2,1,0}=0, \\ 
    2B_4D_4+a_{1,0,0}+c_{0,0,1}=0, \\ 
    2A_3B_4+2A_5D_4+a_{1,0,1}+c_{0,0,2}=0, \\ 
    2A_3A_5-2kA_3'A_5'+a_{1,0,2}+c_{0,0,3}=0, \\ 
    2A_1D_4+2A_5B_4+a_{2,0,0}+c_{1,0,1}=0, \\ 
    2A_1A_3+A_5^2-2kA_1'A_3'-k{A_5'}^2+a_{2,0,1}+c_{1,0,2}=0, \\ 
    2A_1A_5-2kA_1'A_5'+a_{3,0,0}+c_{2,0,1}=0, \\ 
    6A_1B_4-4B_1B_4-c_{2,0,0}=0, \\ 
    -4A_1^2+6A_1B_1+4k{A_1'}^2-2B_1^2-c_{3,0,0}=0, \\ 
    -6A_1B_4+2B_1B_4-c_{2,0,0}=0, \\
    6A_1^2-6A_1B_1-6k{A_1'}^2+B_1^2=0, \\ 
    4A_1^2+2A_1B_1+4k{A_1'}^2-c_{3,0,0}=0, \\ 
    2A_2B_4+2A_4C_4-2B_2C_4-2B_4C_2-c_{0,2,0}=0, \\ 
\end{eqnarray*}

\begin{eqnarray*}
    -2A_2A_4+2A_2B_2+2A_4C_2+2kA_2'A_4'-B_2C_2-c_{0,3,0}=0, \\ 
    -2A_1C_4-2A_4B_4+2B_2B_4-b_{2,0,0}=0, \\ 
    2A_1A_2-2A_1C_2+A_4^2-2A_4B_2-2kA_1'A_2'-k{A_4'}^2+B_2^2=0, \\ 
    -2A_1A_4+2A_1B_2+2kA_1'A_4'-b_{3,0,0}=0, \\ 
    -2A_1A_5+2A_1B_3+2kA_1'A_5'-a_{3,0,0}=0, \\ 
    2A_1A_3-2A_1D_3+A_5^2-2A_5B_3-2kA_1'A_3'-k{A_5'}^2+B_3^2=0, \\ 
    -2A_1D_4-2A_5B_4+2B_3B_4-a_{2,0,0}=0, \\ 
    -2A_3A_5+2A_3B_3+2A_5D_3+2kA_3'A_5'-2B_3D_3-c_{3,0,0}=0, \\ 
    2A_3B_4+2A_5D_4-2B_3D_4-2B_4D_3-c_{0,0,2}=0, \\ 
    2A_2A_6-2kA_2'A_6'+a_{0,3,0}+b_{0,2,1}=0, \\ 
    2A_2A_3+A_6^2-2k{A_2'}{A_3'}-k{A_6'}^2+a_{0,2,1}+b_{0,1,2}=0, \\ 
    2A_2D_4+2A_6C_4+a_{0,2,0}+b_{0,1,1}=0, \\ 
    2A_3A_6-2kA_3'A_6'+a_{0,1,2}+b_{0,0,3}=0, \\ 
    2A_3C_4+2A_6D_4+a_{0,1,1}+b_{0,0,2}=0, \\ 
    2C_4D_4+a_{0,1,0}+b_{0,0,1}=0, \\ 
    -2A_2A_4+2A_2C_1+2kA_2'A_4'-c_{0,3,0}=0, \\ 
     2A_1A_2-2A_2B_1+A_4^2-2A_4C_1-2kA_1'A_2'-k{A_4'}^2+C_1^2=0, \\ 
    -2A_2B_4-2A_4C_4+2C_1C_4-c_{0,2,0}=0, \\ 
    -2A_1A_4+2A_1C_1+2A_4B_1+2kA_1'A_4'-2B_1C_1-b_{3,0,0}=0, \\ 
    2A_1C_4+2A_4B_4-2B_1C_4-2B_4C_1-b_{2,0,0}=0, \\ 
    -4A_2^2+2A_2C_2+4k{A_2'}^2-b_{0,3,0}=0, \\ 
    6A_2^2-6A_2C_2-6k{A_2'}^2+C_2^2=0, \\ 
    -6A_2C_4+2C_2C_4-b_{0,2,0}=0, \\ 
    -4A_2^2+6A_2C_2+4k{A_2'}^2-2C_2^2-b_{0,3,0}=0, \\ 
    6A_2C_4-4C_2C_4-b_{0,2,0}=0, \\ 
    -2A_2A_6+2A_2C_3+2kA_2'A_6'-a_{0,3,0}=0, \\ 
    2A_2A_3-2A_2D_3+A_6^2-2A_6C_3-2kA_2'A_3'-k{A_6'}^2+C_3^2=0, \\ 
    -2A_2D_4-2A_6C_4+2C_3C_4-a_{0,2,0}=0, \\ 
    -2A_3A_6+2A_3C_3+2A_6D_3+2kA_3'A_6'-2C_3D_3-b_{0,0,3}=0, \\ 
    2A_3C_4+2A_6D_4-2C_3D_4-2C_4D_3-b_{0,0,2}=0, \\ 
    -2A_3A_5+2A_3D_1+2kA_3'A_5'-c_{0,0,3}=0, \\ 
    2A_1A_3-2A_3B_1+A_5^2-2A_5D_1-2kA_1'A_3'-k{A_5'}^2+D_1^2=0, \\ 
    -2A_3B_4-2A_5D_4+2D_1D_4-c_{0,0,2}=0, \\ 
    -2A_1A_5+2A_1D_1+2A_5D_1+2kA_1'A_5'-2B_1D_1-a_{3,0,0}=0, \\ 
\end{eqnarray*}

\begin{eqnarray*}
2A_1D_4+2A_5B_4-2B_1D_4-2B_4D_1-a_{2,0,0}=0, \\ 
    -2A_3A_6+2A_3D_2+2kA_3'A_6'-b_{0,0,3}=0, \\ 
    2A_2A_3-2A_3C_2+A_6^2-2A_6D_2-2kA_2'A_3'-k{A_6'}^2+D_2^2=0, \\ 
    -2A_3C_4-2A_6D_4+2D_2D_4-b_{0,0,2}=0, \\ 
    -2A_2A_6+2A_2D_2+2A_6C_2+2kA_2'A_6'-2C_2D_2-a_{0,3,0}=0, \\ 
    2A_2D_4+2A_6C_4-2C_2D_4-2C_4D_2-a_{0,2,0}=0, \\ 
    -4A_3^2+2A_3D_3+4k{A_3'}^2-a_{0,0,3}=0, \\ 
    6A_3^2-6A_3D_3-6k{A_3'}^2+D_3^2=0, \\ 
    -6A_3D_4+2D_3D_4-a_{0,0,2}=0, \\ 
    -4A_3^2+6A_3D_3+4k{A_3'}^2-2D_3^2-a_{0,0,3}=0, \\ 
    6A_3D_4-4D_3D_4-a_{0,0,2}=0, \\ 
    2A_1A_4 - 2A_1B_2 - 2A_1C_1 - 2A_4B_1 -
    2kA_1'A_4' + 2B_1B_2 + 2B_1C_1 + b_{3,0,0} +
    c_{2,1,0}=0, \\ 
    2A_1A_2 - 2A_1C_2 - 2A_2B_1 +
    A_4^2 - 2A_4B_2 - 2A_4C_1 -
    2kA_1'A_2'- \\ - k{A_4'}^2 + 2B_1C_2 +
    B_2^2 + 2B_2C_1 + C_1^2 + b_{2,1,0} +
    c_{1,2,0}=0, \\ 
    - 2A_1C_4 - 2A_4B_4 + 2B_1C_4 +
    2B_2B_4 + 2B_4C_1 + b_{2,0,0} + c_{1,1,0}=0, \\ 
    2A_2A_4 - 2A_2B_2 - 2A_2C_1 - 2A_4C_2 -
    2kA_2'A_4' + 2B_2C_2 + 2C_1C_2 + b_{1,2,0} +
    c_{0,3,0}=0, \\ 
     - 2A_2B_4 - 2A_4C_4 + 2B_2C_4 +
    2B_4C_2 + 2C_1C_4 + b_{1,1,0} + c_{0,2,0}=0, \\ 
    2A_1A_5 - 2A_1B_3 - 2A_1D_1 - 2A_5B_1 -
    2kA_1'A_5' + 2B_1B_3 + 2B_1D_1 + a_{3,0,0} +
    c_{2,0,1}=0, \\ 
    2A_1A_3 - 2A_1D_3 - 2A_3B_1 +
    A_5^2 - 2A_5B_3 - 2A_5D_1 -
    2kA_1'A_3'- \\ - k{A_5'}^2 + 2B_1D_3 +
    B_3^2 + 2B_3D_1 + D_1^2 + a_{2,0,1} +
    c_{1,0,2}=0, \\ 
    2A_1D_4 - 2A_5B_4 + 2B_1D_4 +
    2B_3B_4 + 2B_4D_1 + a_{2,0,0} + c_{1,0,1}=0, \\ 
    2A_3A_5 - 2A_3B_3 - 2A_3D_1 - 2A_5D_3 -
    2kA_5' + 2B_3D_3 + 2D_1D_3 + a_{1,0,2} +
    c_{0,0,3}=0, \\ 
    - 2A_3B_4 - 2A_5D_4 + 2B_3D_4 +2B_4D_3 + 2D_1D_4 + a_{1,0,1} + c_{0,0,2}=0, \\ 
    2A_2A_6 - 2A_2C_3 - 2A_2D_2 - 2A_6C_2 -
    2kA_2'A_6' + 2C_2C_3 + 2C_2D_2 + a_{0,3,0} +
    b_{0,2,1}=0, \\ 
    2A_2A_3 - 2A_2D_3 - 2A_3C_2 +
    A_6^2 - 2A_6C_3 - 2A_6D_2 -
    2kA_2'A_3'- \\ - k{A_6'}^2 + 2C_2D_3 +
    C3^2 + 2C_3D_2 + D_2^2 + a_{0,2,1} +
    b_{0,1,2}=0, \\ 
   - 2A_2D_4 - 2A_6C_4 +2C_2D_4 +
    2C_3C_4 + 2C_4D_2 + a_{0,2,0} + b_{0,1,1}=0, \\ 
    2A_3A_6 - 2A_3C_3 - 2A_3D_2 - 2A_6D_3 -
    2kA_3'A_6' + 2C_3D_3 + 2D_2D_3 + a_{0,1,2} +
    b_{0,0,3}=0, \\ 
    - 2A_3C_4 - 2A_6D_4 + 2C_3D_4 +
    2C_4D_3 + 2D_2D_4 + a_{0,1,1} + b_{0,0,2}=0, \\ 
    2A_1A_6 + 2A_4A_5 - 2kA_1'A_6 -
    2kA_4'A_5' + a_{2,1,0} + b_{2,0,1} + c_{1,1,1}=0, \\ 
    2A_2A_5 + 2A_4A_6 -
    2kA_2'A_5' - 2kA_4'A_6' + a_{1,2,0} +
    b_{1,1,1} + c_{0,2,1}=0, \\ 
    2A_3A_4 + 2A_5A_6
    - 2kA_3'A_4' - 2kA_5'A_6' + a_{1,1,1} +
    b_{1,0,2} + c_{0,1,2}=0, \\ 
    2A_4D_4 + 2A_5C_4 +
    2A_6B_4 + a_{1,1,0} + b_{1,0,1} + c_{0,1,1}=0, \\ 
    \end{eqnarray*}

    \begin{eqnarray*}
    6A_1A_4 + 2A_1C_1 + 2A_4B_1 +
    6kA_1'A_4' - c_{2,1,0}=0, \\ 
    4A_1A_2 + 2A_2B_1 - 2A_4^2 +
    2A_4C_1 + 4kA_1'A_2' + 2k{A_4'}^2 -
    c_{1,2,0}=0, \\ 
    6A_1A_4 - 4A_1C_1 -
    4A_4B_1 - 6kA_1'A_4' + 2B_1C_1=0, \\ 
    -4A_1C_4 - 4A_4B_4 + 2B_1C_4 + 2B_4C_1 -
    c_{1,1,0}=0, \\ 
    -4A_1A_2 + 2A_1C_2 - 2A_4^2 +
    2A_4B_2 + 4kA_1'A_2' + 2k{A_4'}^2 -
   b_{2,1,0}=0, \\ 
    6A_2A_4 +
    2A_2B_2 + 2A_4C_2 + 6kA_2'A_4' -
    b_{1,2,0}=0, \\ 
    6A_2A_4 - 4A_2B_2 -
    4A_4C_2 - 6kA_2'A_4' + 2B_2C_2=0, \\ 
    -4A_2B_4 - 4A_4C_4 + 2B_2C_4 + 2B_4C_2 -
    b_{1,1,0}=0, \\ 
    -2A_1A_6 + 2A_1C_3 - 2A_4A_5 +
    2A_4B_3 + 2kA_1'A_6' + 2kA_4'A_5' -
    a_{2,1,0}=0, \\ 
    2A_2A_5 +
    2A_2B_3 - 2A_4A_6 + 2A_4C_3 +
    2kA_2'A_5' + 2kA_4'A_6' - a_{1,2,0}=0, \\ 
    2A_3A_4 - 2A_4D_3 + 2A_5A_6 -
    2A_5C_3 - 2A_6B_3 - 2kA_3'A_4' -
    2kA_5'A_6' + 2B_3C_3=0, \\ 
    - 2A_4D_4 -
    2A_5C_4 - 2A_6B_4 + 2B_3C_4 + 2B_4C_3 -
    a_{1,1,0}=0, \\ 
    -6A_1A_5 + 2A_1D_1 + 2A_5B_1 +
    6kA_1'A_5' - c_{2,0,1}=0, \\ 
    4A_1A_3 + 2A_3B_1 - 2A5^2 +
    2A_5D_1 + 4kA_1'A_3' + 2k{A_5'}^2 -
    c_{1,0,2}=0, \\ 
    6A_1A_5 - 4A_1D_1 -
    4A_5B_1 - 6kA_1'A_5' + 2B_1D_1=0, \\ 
    4A_1D_4 - 4A_5B_4 + 2B_1D_4 + 2B_4D_1 -
    c_{1,0,1}=0, \\ 
    -2A_1A_6 + 2A_1D_2 - 2A_4A_5 +
    2A_5B_2 + 2kA_1'A_6' + 2kA_4'A_5' -
    b_{2,0,1}=0, \\ 
    2A_3B_2 - 2A_5A_6 + 2A_5D_2 +
    2kA_3'A_4' + 2kA_5'A_6' - b_{1,0,2}=0, \\ 
    2A_2A_5 + 2A_4A_6 - 2A_4D_2 -
    2A_5C_2 - 2A_6B_2 - 2kA_2'A_5' -
    2kA_4'A_6' + 2B_2D_2=0, \\ 
    - 2A_4D_4 -2A_5C_4 - 2A_6B_4 + 2B_2D_4 + 2B_4D_2 -b_{1,0,1}=0, \\ 
    -4A_1A_3 + 2A_1D_3 - 2A_5^2 +
    2A_5B_3 + 4kA_1'A_3' + 2k{A_5'}^2 -
    a_{2,0,1}=0, \\ 
    6A_3A_5 +
    2A_3B_3 + 2A_5D_3 + 6kA_3'A_5' -
    a_{1,0,2}=0, \\ 
    6A_3A_5 - 4A_3B_3 -
    4A_5D_3 - 6kA_3'A_5' + 2B_3D_3=0, \\ 
    -4A_3B_4 - 4A_5D_4 + 2B_3D_4 + 2B_4D_3 -
    a_{1,0,1}=0, \\ 
    6A_1A_4 - 6A_1B_2 - 2A_1C_1 -
    2A_4B_1 - 6kA_1'A_4' + 2B_1B_2=0, \\ 
    6A_1A_4 + 6A_1B_2 +
    4A_1C_1 + 4A_4B_1 + 6kA_1'A_4' -
    4B_1B_2 - 2B_1C_1 - c_{2,1,0}=0, \\ 
    -4A_1A_2 + 4A_1C_2 + 2A_2B_1 -
    2A_4^2 + 4A_4B_2 + 2A_4C_1 + \\
    4kA_1'A_2' + 2k{A_4'}^2 - 2B_1C_2 -
    2B_2^2 - 2B_2C_1 - c_{1,2,0}=0, \\ 
    4A_1C_4
    + 4A_4B_4 - 2B_1C_4 - 4B_2B_4 - 2B_4C_1
    - c_{1,1,0}=0, \\ 
    6A_1A_5 - 6A_1B_3 - 2A_1D_1 -
    2A_5B_1 - 6kA_1'A_5' + 2B_1B_3=0, \\ 
    6A_1A_5 + 6A_1B_3 +
    4A_1D_1 + 4A_5B_1 + 6kA_1'A_5' -
    4B_1B_3 - 2B_1D_1 - c_{2,0,1}=0, \\ 
    -4A_1A_3 + 4A_1D_3 + 2A_3B_1 -
    2A_5^2 + 4A_5B_3 + 2A_5D_1 + \\
    4kA_1'A_3' + 2k{A_5'}^2 - 2B_1D_3 -
    2B_3^2 - 2B_3D_1 - c_{1,0,2}=0, \\ 
    \end{eqnarray*}

    \begin{eqnarray*} 
    4A_1D_4
    + 4A_5B_4 - 2B_1D_4 - 4B_3B_4 - 2B_4D_1
    - c_{1,0,1}=0, \\ 
    2A_1A_6 - 2A_1C_3 - 2A_1D_2 +
    2A_4A_5 - 2A_4B_3 - 2A_5B_2 -
    2kA_1'A_6' - 2kA_4'A_5' + 2B_2B_3=0, \\ 
    2A_2A_5 + 2A_2B_3 -
    2A_4A_6 + 2A_4C_3 + 2A_4D_2 +
    2A_5C_2 + 2A_6B_2 + \\ 2kA_2'A_5' +
    2kA_4'A_6' - 2B_2C_3 - 2B_2D_2 -
    2B_3C_2 - c_{0,2,1}=0, \\ 
    - 2A_3A_4 +
    2A_3B_2 + 2A_4D_3 - 2A_5A_6 +
    2A_5C_3 + 2A_5D_2 + 2A_6B_3 +\\
    2kA_3'A_4' + 2kA_5'A_6' - 2B_2D_3 -
    2B_3C_3 - 2B_3D_2 - c_{0,1,2}=0, \\ 
    2A_4D_4
    + 2A_5C_4 + 2A_6B_4 - 2B_2D_4 - 2B_3C_4
    - 2B_4C_3 - 2B_4D_2 - c_{0,1,1}=0, \\ 
    2A_2A_5 + 2A_2D_1 - 2A_4A_6 +
    2A_6C_1 + 2kA_2'A_5' + 2kA_4'A_6' -
    c_{0,2,1}=0, \\ 
    2A_3A_4 +
    2A_3C_1 - 2A_5A_6 + 2A_6D_1 +
    2kA_3'A_4' + 2kA_5'A_6' - c_{0,1,2}=0, \\ 
    2A_1A_6 + 2A_4A_5 - 2A_4D_1 -
    2A_5C_1 - 2A_6B_1 - 2kA_1'A_6' -
    2kA_4'A_5' + 2C_1D_1=0, \\ 
    - 2A_4D_4 -
    2A_5C_4 - 2A_6B_4 + 2C_1D_4 + 2C_4D_1 -
    c_{0,1,1}=0, \\ 
    -6A_2A_6 + 2A_2D_2 + 2A_6C_2 +
    6kA_2'A_6' - b_{0,2,1}=0, \\ 
    4A_2A_3 + 2A_3C_2 - 2A_6^2 +
    2A_6D_2 + 4kA_2'A_3' + 2k{A_6'}^2 -
    b_{0,1,2}=0, \\ 
    6A_2A_6 - 4A_2D_2 -
    4A_6C_2 - 6kA_2'A_6' + 2C_2D_2=0, \\ 
    -4A_2D_4 - 4A_6C_4 + 2C_2D_4 + 2C_4D_2 -
    b_{0,1,1}=0, \\ 
    -4A_2A_3 + 2A_2D_3 - 2A_6^2 +
    2A_6C_3 + 4kA_2'A_3' + 2k{A_6'}^2 -
    a_{0,2,1}=0, \\ 
    - 6A_3A_6 +
    2A_3C_3 + 2A_6D_3 + 6kA_3'A_6' -
    a_{0,1,2}=0, \\ 
     6A_3A_6 - 4A_3C_3 -
    4A_6D_3 - 6kA_3'A_6' + 2C_3D_3=0, \\ 
    -4A_3C_4 - 4A_6D_4 + 2C_3D_4 + 2C_4D_3 -
    a_{0,1,1}=0, \\ 
    6A_2A_4 - 2A_2B_2 - 6A_2C_1 -
    2A_4C_2 - 6kA_2'A_4' + 2C_1C_2=0, \\ 
    -4A_1A_2 + 2A_1C_2 +
    4A_2B_1 - 2A_4^2 + 2A_4B_2 +
    4A_4C_1 + \\ 4kA_1'A_2' + 2k{A_4'}^2 -
    2B_1C_2 - 2B_2C_1 - 2C_1^2 -
    b_{2,1,0}=0, \\ 
    - 6A_2A_4 + 4A_2B_2 +
    6A_2C_1 + 4A_4C_2 + 6kA_2'A_4' -
    2B_2C_2 - 4C_1C_2 - b_{1,2,0}=0, \\ 
    4A_2B_4
    + 4A_4C_4 - 2B_2C_4 - 2B_4C_2 - 4C_1C_4
    - b_{1,1,0}=0, \\ 
    2A_2A_5 - 2A_2B_3 - 2A_2D_1 +
    2A_4A_6 - 2A_4C_3 - 2A_6C_1 -
    2kA_2'A_5' - 2kA_4'A_6' + 2C_1C_3=0, \\ 
    - 2A_1A_6 + 2A_1C_3 -
    2A_4A_5 + 2A_4B_3 + 2A_4D_1 +
    2A_5C_1 + \\ 2A_6B_1 + 2kA_1'A_6' +
    2kA_4'A_5' - 2B_1C_3 - 2B_3C_1 -
    2C_1D_1 - b_{2,0,1}=0, \\ 
    - 2A_3A_4 +
    2A_3C_1 + 2A_4D_3 - 2A_5A_6 +
    2A_5C_3 + 2A_6B_3 + 2A_6D_1 +\\
    2kA_3'A_4' + 2kA_5'A_6' - 2B_3C_3 -
    2C_1D_3 - 2C_3D_1 - b_{1,0,2}=0, \\ 
    2A_4D_4
    + 2A_5C_4 + 2A_6B_4 - 2B_3C_4 - 2B_4C_3
    - 2C_1D_4 - 2C_4D_1 - b_{1,0,1}=0, \\ 
    6A_2A_6 - 6A_2C_3 - 2A_2D_2 -
    2A_6C_2 - 6kA_2'A_6' + 2C_2C_3=0, \\ 
    - 6A_2A_6 + 6A_2C_3 +
    4A_2D_2 + 4A_6C_2 + 6kA_2'A_6' -
    4C_2C_3 - 2C_2D_2 - b_{0,2,1}=0, \\ 
    \end{eqnarray*}

    \begin{eqnarray*}
    -4A_2A_3 + 4A_2D_3 + 2A_3C_2 -
    2A_6^2 + 4A_6C_3 + 2A_6D_2 +\\
    4kA_2'A_3' + 2k{A_6'}^2 - 2C_2D_3 -
    2C_3^2 - 2C_3D_2 - b_{0,1,2}=0, \\ 
    4A_2D_4
    + 4A_6C_4 - 2C_2D_4 - 4C_3C_4 - 2C_4D_2
    - b_{0,1,1}=0, \\ 
    2A_3A_4 - 2A_3B_2 - 2A_3C_1 +
    2A_5A_6 - 2A_5D_2 - \\ 2A_6D_1 -
    2kA_3'A_4' - 2kA_5'A_6' + 2D_1D_2=0, \\ 
    - 2A_1A_6 + 2A_1D_2 -
    2A_4A_5 + 2A_4D_1 + 2A_5B_2 +
    2A_5C_1 + 2A_6B_1 + \\ 2kA_1'A_6' +
    2kA_4'A_5' - 2B_1D_2 - 2B_2D_1 -
    2C_1D_1 - a_{2,1,0}=0, \\ 
    - 2A_2A_5 +
    2A_2D_1 - 2A_4A_6 + 2A_4D_2 +
    2A_5C_2 + 2A_6B_2 + 2A_6C_1 + \\
    2kA_2'A_5' + 2kA_4'A_6' - 2B_2D_2 -
    2C_1D_2 - 2C_2D_1 - a_{1,2,0}=0, \\ 
    2A_4D_4
    + 2A_5C_4 + 2A_6B_4 - 2B_2D_4 - 2B_4D_2
    - 2C_1D_4 - 2C_4D_1 - a_{1,1,0}=0, \\ 
    6A_3A_5 - 2A_3B_3 - 6A_3D_1 -
    2A_5D_3 - 6kA_3'A_5' + 2D_1D_3=0, \\ 
    - 4A_1A_3 + 2A_1D_3 +
    4A_3B_1 - 2A_5^2 + 2A_5B_3 +
    4A_5D_1 + \\4kA_1'A_3' + 2k{A_5'}^2 -
    2B_1D_3 - 2B_3D_1 - 2D_1^2 -
    a_{2,0,1}=0, \\ 
    - 6A_3A_5 + 4A_3B_3 +
    6A_3D_1 + 4A_5D_3 + 6kA_3'A_5' -
    2B_3D_3 - 4D_1D_3 - a_{1,0,2}=0, \\ 
     4A_3B_4
    + 4A_5D_4 - 2B_3D_4 - 2B_4D_3 - 4D_1D_4
    - a_{1,0,1}=0, \\ 
    6A_3A_6 - 2A_3C_3 - 6A_3D_2 -
    2A_6D_3 - 6kA_3'A_6' + 2D_2D_3=0, \\ 
    - 4A_2A_3 + 2A_2D_3 +
    4A_3C_2 - 2A_6^2 + 2A_6C_3 +
    4A_6D_2 + \\ 4kA_2'A_3' + 2k{A_6'}^2-
    2C_2D_3 - 2C_3D_2 - 2D_2^2 -
    a_{0,2,1}=0, \\ 
    - 6A_3A_6 + 4A_3C_3 +
    6A_3D_2 + 4A_6D_3 + 6kA_3'A_6' -
    2C_3D_3 - 4D_2D_3 - a_{0,1,2}=0, \\ 
    4A_3C_4
    + 4A_6D_4 - 2C_3D_4 - 2C_4D_3 - 4D_2D_4
    - a_{0,1,1}=0, \\ 
    2A_1A_6 -
    2A_1C_3 - 2A_1D_2 + 2A_4A_5 -
    2A_4B_3 - 2A_4D_1 -  2A_5B_2 -\\
    2A_5C_1 - 2A_6B_1 - 2kA_1'A_6'-
    2kA_4'A_5' + 2B_1C_3 + 2B_1D_2 +
    2B_2B_3 + \\2B_2D_1 + 2B_3C_1 +
    2C_1D_1 + a_{2,1,0} + b_{2,0,1} + c_{1,1,1}=0, \\ 
    2A_2A_5 - 2A_2B_3 -
    2A_2D_1 + 2A_4A_6 - 2A_4C_3 -
    2A_4D_2 - 2A_5C_2 - \\ 2A_6B_2 -
    2A_6C_1 - 2kA_2'A_5' - 2kA_4'A_6' +
    2B_2C_3 + 2B_2D_2 + 2B_3C_2 +\\
    2C_1C_3 + 2C_1D_2 + 2C_2D_1 +
    a_{1,2,0} + b_{1,1,1} +c_{0,2,1}=0, \\ 
    2A_3A_4 -
    2A_3B_2 - 2A_3C_1 - 2A_4D_3 +
    2A_5A_6 - 2A_5C_3 -\\ 2A_5D_2 -
    2A_6B_3 - 2A_6D_1 - 2kA_3'A_4' -
    2A_5'A_6 + 2B_2D_3 +\\ 2B_3C_3 +
    2B_3D_2 + 2C_1D_3 + 2C_3D_1 +
    2D_1D_2 + a_{1,1,1} + b_{1,0,2} + c_{0,1,2}=0, \\ 
    -2A_4D_4 - 2A_5C_4 - 2A_6B_4 + 2B_2D_4 +
    2B_3C_4 + 2B_4C_3 + \\ 2B_4D_2 + 2C_1D_4 +
    2C_4D_1 + a_{1,1,0} + b_{1,0,1} + c_{0,1,1}=0, \\ 
    - 4A_1A_6 -
    4A_4A_5 + 2A_4D_1 + 2A_5C_1 +
    2A_6B_1 + 4kA_1'A_6' + 4kA_4'A_5' -
    c_{1,1,1}=0, \\ 
    - 4A_2A_5 -
    4A_4A_6 + 2A_4D_2 + 2A_5C_2 +
    2A_6B_2 + 4kA_2'A_5' + 4kA_4'A_6' -
    b_{1,1,1}=0, \\ 
    - 4A_3A_4 +
    2A_4D_3 - 4A_5A_6 + 2A_5C_3 +
    2A_6B_3 + 4kA_3'A_4' + 4A_5'A_6' -
    a_{1,1,1}=0, \\ 
    \end{eqnarray*}

    \begin{eqnarray*}
    4A_1A_6 -
    4A_1C_3 + 4A_4A_5 - 4A_4B_3 -
    2A_4D_1 - 2A_5C_1 - \\ 2A_6B_1 -
    4kA_1'A_6' - 4kA_4'A_5' + 2B_1C_3 +
    2B_3C_1=0, \\ 
    8A_1A_2 - 4A_1C_2 - 4A_2B_1 +
    4A_4^2 - 4A_4B_2 - \\4A_4C_1 -
    8kA_1'A_2' - 4k{A_4'}^2 + 2B_1C_2 +
    2B_2C_1=0, \\ 
    4A_2A_5 -
    4A_2B_3 + 4A_4A_6 - 4A_4C_3 -
    2A_4D_2 - 2A_5C_2 -\\ 2A_6B_2 -
    4kA_2'A_5' - 4kA_4'A_6' + 2B_2C_3 +
    2B_3C_2=0, \\ 
    4A_1A_6 -
    4A_1D_2 + 4A_4A_5 - 2A_4D_1 -
    4A_5B_2 - \\2A_5C_1 - 2A_6B_1 -
    4kA_1'A_6' - 4kA_4'A_5' + 2B_1D_2 +
    2B_2D_1=0, \\ 
    8A_1A_3 - 4A_1D_3 - 4A_3B_1 +
    4A_5^2 - 4A_5B_3 - \\4A_5D_1 -
    8kA_1'A_3' - 4k{A_5'}^2 + 2B_1D_3 +
    2B_3D_1=0, \\ 
    4A_3A_4 -
    4A_3B_2 - 2A_4D_3 + 4A_5A_6 -
    2A_5C_3 - 4A_5D_2 - \\ 2A_6B_3 -
    4kA_3'A_4' - 4kA_5'A_6' + 2B_2D_3 +
    2B_3D_2=0, \\ 
    -4A_1A_6 + 4A_1C_3 + 4A_1D_2 -
    4A_4A_5 + 4A_4B_3 + 2A_4D_1 +
    4A_5B_2 + 2A_5C_1 + \\ 2A_6B_1 +
    4kA_1'A_6' + 4kA_4'A_5' - 2B_1C_3 -
    2B_1D_2 - 4B_2B_3 - 2B_2D_1 -
    2B_3C_1 - c_{1,1,1}=0, \\ 
    4A_2A_5 -
    4A_2D_1 + 4A_4A_6 - 2A_4D_2 -
    2A_5C_2 - \\ 2A_6B_2 - 4A_6C_1 -
    4kA_2'A_5' - 4kA_4'A_6'+ 2C_1D_2 +
    2C_2D_1=0, \\ 
    4A_3A_4 -
    4A_3C_1 - 2A_4D_3 + 4A_5A_6 -
    2A_5C_3 - \\ 2A_6B_3 - 4A_6D_1 -
    4kA_3'A_4' - 4kA_5'A_6' + 2C_1D_3 +
    2C_3D_1=0, \\ 
    8A_2A_3 - 4A_2D_3 - 4A_3C_2 +
    4A_6^2 - 4A_6C_3 - \\ 4A_6D_2 -
    8kA_2'A_3' - 4k{A_6'}^2 + 2C_2D_3 +
    2C_3D_2=0, \\ 
    - 4A_2A_5 + 4A_2B_3 +
    4A_2D_1 - 4A_4A_6 + 4A_4C_3 +
    2A_4D_2 + 2A_5C_2 + 2A_6B_2 +
    4A_6C_1 + \\ 4kA_2'A_5' + 4kA_4'A_6'-
    2B_2C_3 - 2B_3C_2 - 4C_1C_3 -
    2C_1D_2 - 2C_2D_1 - b_{1,1,1}=0, \\ 
    -4A_3A_4 + 4A_3B_2 + 4A_3C_1 +
    2A_4D_3 - 4A_5A_6 + 2A_5C_3 +
    4A_5D_2 + 2A_6B_3 + 4A_6D_1 +\\
    4kA_3'A_4' + 4kA_5'A_6' - 2B_2D_3 -
    2B_3D_2 - 2C_1D_3 - 2C_3D_1 -
    4D_1D_2 - a_{1,1,1}=0. \\ 
    \end{eqnarray*}

   \textbf{CASE $p=2$}

    \begin{eqnarray*}
    A_1^2 + A_1A_1' + {A_1'}^2(\xi^2 + \xi) + c_{3,0,0}=0, \\ 
    A_1A_4' +A_4A_1' + b_{3,0,0} + c_{2,1,0}=0, \\
    A_1A_5' + A_5A_1' + a_{3,0,0}+
        c_{2,0,1}=0, \\
        A_1'B_1 +c_{3,0,0}=0,\\
        A_1A_4' + A_4A_1' + A_1'B_2 +
        b_{3,0,0}=0, \\
        A_1A_5' + A_5A_1' + A_1'B_3 + a_{3,0,0}=0, \\
        A_1'B_4 + c_{2,0,0}=0, \\
     A_1A_2' + A_2A_1' + A_4^2 + A_4A_4' + {A_4'}^2(\xi^2 +
        \xi) +b_{2,1,0} + c_{1,2,0}=0,\\
        A_1A_6' +A_4A_5' +
        A_5A_4' + A_6A_1' + a_{2,1,0} +b_{2,0,1} + c_{1,1,1}=0, \\
        A_1A_4' +A_4A_1' + A_1'C_1 + A_4'B_1 + c_{2,1,0}=0,\\
        A_1'C_2 + A_4'B_2 + b_{2,1,0}=0,\\
        A_1A_6' + A_4A_5' +
        A_5A_4' + A_6A_1' + A_1'C_3 + A_4'B_3 + a_{2,1,0}=0,\\
        A_1'C_4 + A_4'B_4 + b_{2,0,0} + c_{1,1,0}=0,\\ 
        A_1A_3' + A_3A_1' +
        A_5^2 + A_5A_5' + {A_5'}^2(\xi^2 + \xi) + a_{2,0,1} +c_{1,0,2}=0,\\ 
        A_1A_5' + A_5A_1' + A_1'D_1 + A_5'B_1 +
        c_{2,0,1}=0,\\ 
        A_1A_6' + A_4A_5' + A_5A_4' + A_6A_1' +
        A_1'D_2 + A_5'B_2 + b_{2,0,1}=0,\\ 
        A_1'D_3 + A_5'B_3 +
        a_{2,0,1}=0,\\ 
        A_1'D_4 + A_5'B_4 + a_{2,0,0} + c_{1,0,1}=0,\\ 
        A_1'B_1 +
        B_1^2=0,\\ 
        A_1A_4' + A_4A_1' + A_1'B_2 + A_1'C_1 +
       A_4'B_1=0,\\ 
       A_1A_5' + A_5A_1' + A_1'B_3 + A_1'D_1
        + A_5'B_1=0,\\ 
        A_1A_2' + A_2A_1' +
        A_4^2 + A_4A_4' + A_1'C_2 + {A_4'}^2(\xi^2 + \xi)
        + A_4'B_2 + B_2^2=0,\\ 
         A_1A_6' + A_4A_5' + A_5A_4' +
        A_6A_1' + A_1'C_3 + A_1'D_2 + A_4'B_3 + A_5'B_2=0,\\
        A_1'C_4 + A_4'B_4 + b_{2,0,0}=0,\\
        A_1A_3' + A_3A_1' +
        A_5^2 + A_5A_5' + A_1'D_3 + {A_5'}^2(\xi^2 + \xi)
        + A_5'B_3 + B_3^2=0,\\ 
        A_1'D_4 + A_5'B_4 + a_{2,0,0}=0,\\
        B_4^2 +c_{1,0,0}=0, \\ 
          A_2A_4' + A_4A_2' + b_{1,2,0} + c_{0,3,0}=0,\\
          A_2A_5' +
        A_4A_6' + A_5A_2' + A_6A_4' + a_{1,2,0} +
        b_{1,1,1} + c_{0,2,1}=0, \\ 
         A_2'B_1 + A_4'C_1 +
        c_{1,2,0}=0, \\
        A_2A_4' + A_4A_2' + A_2'B_2 +
        A_4'C_2 + b_{1,2,0}=0, \\
         A_2A_5' + A_4A_6' +
        A_5A_2' + A_6A_4' + A_2'B_3 + A_4'C_3 +
        a_{1,2,0}=0, \\
        A_2'B_4 + A_4'C_4 + b_{1,1,0} + c_{0,2,0}=0,\\
        A_3A_4' + A_4A_3' + A_5A_6'+ A_6A_5' +
       a_{1,1,1} + b_{1,0,2} + c_{0,1,2}=0, \\
       \end{eqnarray*}
        \begin{eqnarray*}
        A_4'D_1 +
        A_5'C_1 + A_6'B_1 + c_{1,1,1}=0, \\ 
        A_4'D_2 +
        A_5'C_2 + A_6'B_2 + b_{1,1,1}=0, \\
       A_4'D_3 +
        A_5'C_3 + A_6'B_3 + a_{1,1,1}=0, \\
        A_4'D_4 +
        A_5'C_4 + A_6'B_4 + a_{1,1,0} + b_{1,0,1} + c_{0,1,1}=0, \\
        A_1A_4' + A_4A_1'=0, \\
        A_4'D_1 + A_5'C_1 +
        A_6'B_1=0, \\
        c_{1,1,0}=0,\\
        A_2A_4' + A_4A_2=0, \\
        A_4'D_2 + A_5'C_2 + A_6'B_2=0, \\
        b_{1,1,0}=0, \\
        A_3A_4' + A_4A_3' + A_5A_6' + A_6A_5' +
        A_4'D_3 + A_5'C_3 + A_6'B_3=0, \\
        A_4'D_4 +
        A_5'C_4 + A_6'B_4 + a_{1,1,0}=0, \\
        b_{1,0,0} + c_{0,1,0}=0, \\
        A_5' + A_5A_3' + a_{1,0,2} + c_{0,0,3}=0, \\
        A_3'B_1 +
        A_5'D_1 + c_{1,0,2}=0, \\
        A_3A_4' +A_4A_3' +
        A_5A_6' + A_6A_5' + A_3'B_2 + A_5'D_2 +
        b_{1,0,2}=0, \\
        A_3A_5' + A_5A_3' + A_3'B_3 +
        A_5'D_3 + a_{1,0,2}=0, \\
        A_3'B_4 +A_5'D_4 + a_{1,0,1} +
        c_{0,0,2}=0,\\
        A_1A_5' + A_5A_1'=0,\\
        A_4'D_1 +
        A_5'C_1 + A_6'B_1=0,\\
         c_{1,0,1}=0, \\
        A_2A_5' +
        A_4A_6' + A_5A_2' + A_6A_4'+ A_4'D_2 +
        A_5'C_2 + A_6'B_2=0,\\
        A_4'D_3 + A_5'C_3 +
       A_6'B_3=0,\\
     A_4'D_4 + A_5C_4 + A_6'B_4 +
        b_{1,0,1}=0,\\
        A_3A_5' + A_5A_3'=0,\\
         a_{1,0,1}=0,\\
        a_{1,0,0} +
        c_{0,0,1}=0,\\
         A_1A_4' + A_4A_1' +
        A_1'B_2 + c_{2,1,0}=0,\\
        A_1A_5' + A_5A_1' +
        A_1'B_3 + c_{2,0,1}=0, \\
         A_2'B_1
        + A_4'C_1 + c_{1,2,0}=0, \\
        A_4'D_1 + A_5'C_1 +
        A_6'B_1 + c_{1,1,1}=0,\\
        A_3'B_1 +
        A_5'D_1 + c_{1,0,2}=0, \\
        A_2A_4' + A_4A_2' +
        A_2'B_2 + A_4'C_2 + c_{0,3,0}=0,\\
        \end{eqnarray*}
        \begin{eqnarray*}
        A_2A_5' + A_4A_6'
        + A_5A_2' + A_6A_4' + A_2'B_3 + A_4'C_3 +
        A_4'D_2 + A_5'C_2 + A_6'B_2 + c_{0,2,1}=0,\\
        A_2'B_4 + A_4'C_4 + c_{0,2,0}=0,\\
        A_3A_4' + A_4A_3'
       A_5A_6' + A_6A_5' + A_3'B_2 + A_4'D_3 +
    A_5'C_3 + A_5'D_2 + A_6'B_3 + c_{0,1,2}=0, \\
        A_4'D_4 + A_5'C_4 + A_6'B_4 + c_{0,1,1}=0, \\
        A_3A_5' +A_3' + A_3'B_3 + A_5'D_3 +
        c_{0,0,3}=0, \\
        A_3'B_4 + A_5'D_4 + c_{0,0,2}=0, \\
        A_2^2 + A_2A_2' + {A_2'}^2(\xi^2 + \xi) + b_{0,3,0}=0, \\
        A_2A_6' + A_6A_2' + a_{0,3,0} + b_{0,2,1}=0, \\
        A_2A_4' + A_4A_2' + A_2'C_1 + c_{0,3,0}=0, \\
        A_2'C_2 + b_{0,3,0}=0, \\  
        A_2A_6' + A_6A_2' +
        A_2'C_3 + a_{0,3,0}=0, \\
        A_2'C_4 + b_{0,2,0}=0, \\
        A_2A_3' + A_3A_2' + A_6^2 + A_6A_6'
        + {A_6'}^2(\xi^2 + \xi) + a_{0,2,1} +
        b_{0,1,2}=0, \\
        A_2A_5' + A_4A_6' +
        A_5A_2'+ A_6A_4' + A_2'D_1 +
        A_6'C_1 + c_{0,2,1}=0, \\ 
        A_2A_6' +
        A_6A_2' + A_2'D_2 + A_6'C_2 +
        b_{0,2,1}=0, \\
        A_2'D_3 + A_6'C_3 +
        a_{0,2,1}=0, \\
        A_2'D_4 + A_6'C_4 + a_{0,2,0} +
        b_{0,1,1}=0, \\
         A_1'A_2' + A_2A_1' + A_4^2 +
        A_4A_4' + A_2'B_1 + {A_4'}^2(\xi^2 +
       \xi) + A_4'C_1 + C_1^2=0, \\
        A_2A_4' + A_4A_2' + A_2'B_2 +
        A_2'C_1 + A_4'C_2=0, \\
        A_2A_5' +
        A_4A_6' + A_5A_2' + A_6A_4' +
        A_2'B_3 + A_2'D_1 + A_4'C_3 +
        A_6'C_1=0, \\
        A_2'B_4 + A_4'C_4 + c_{0,2,0}=0, \\
        A_2'C_2 + C_2^2=0, \\ 
         A_2A_6' +
        A_6A_2' + A_2'C_3 + A_2'D_2 +
        A_6'C_2=0, \\
        A_2A_3' +
        A_3A_2' + A_6^2 + A_6A_6' + A_2'D_3
        + {A_6'}^2(\xi^2 + \xi) + A_6'C_3 +
        C_3^2=0, \\
        A_2'D_4 + A_6'C_4 + a_{0,2,0}=0, \\
        C_4^2 + b_{0,1,0}=0, \\
        A_3A_6' + A_6A_3' + a_{0,1,2} +
        b_{0,0,3}=0, \\
        A_3A_4' + A_4A_3' + A_5A_6'
        + A_6A_5' + A_3'C_1 + A_6'D_1 +
        c_{0,1,2}=0, \\
        A_3'C_2 + A_6'D_2 +
        b_{0,1,2}=0, \\
        A_3A_6' + A_6A_3' +
        A_3'C_3 + A_6'D_3 + a_{0,1,2}=0, \\
        A_3'C_4
        + A_6'D_4 + a_{0,1,1} + b_{0,0,2}=0, \\
        A_1A_6' +
        A_4A_5' + A_5A_4' + A_6A_1' +
        A_4'D_1 + A_5'C_1 + A_6'B_1=0, \\
       A_4'D_2 + A_5'C_2 + A_6'B_2=0, \\
       \end{eqnarray*}
        \begin{eqnarray*}
        A_4'D_3 + A_5'C_3 + A_6'B_3=0, \\
        A_4'D_4 + A_5'C_4 + A_6'B_4 + c_{0,1,1}=0, \\
        A_2A_6' + A_6A_2'=0, \\
        A_3A_6' + A_6A_3'=0, \\ a_{0,1,1}=0, \\
        a_{0,1,0} +
        b_{0,0,1}=0, \\
        A_1A_4' +A_4A_1' + A_1'C_1 +
        A_4'B_1 + b_{3,0,0}=0, \\ 
        A_1'C_2 + A_4'B_2 +
        b_{2,1,0}=0, \\
        A_1A_6' + A_4A_5' +
        A_5A_4' + A_6A_1' + A_1'C_3 +
       A_4'B_3 + A_4'D_1 + A_5'C_1 +
        A_6'B_1 + b_{2,0,1}=0, \\ A_1'C_4 + A_4'B_4 +
       b_{2,0,0}=0, \\
        A_2A_4' + A_4A_2' + A_2'C_1
        + b_{1,2,0}=0, \\
        A_4'D_2 + A_5'C_2 +
        A_6'B_2 + b_{1,1,1}=0, \\ b_{1,1,0}=0, \\ 
        A_3A_4' + A_4A_3' + A_5A_6' +
        A_6A_5' + A_3'C_1 + A_4'D_3 +
        A_5'C_3 + A_6'B_3 + A_6'D_1 +
        b_{1,0,2}=0,\\
         A_4D_4 + A_5'C_4 + A_6'B_4 +
        b_{1,0,1}=0, \\
        b_{1,0,0} + c_{0,1,0}=0, \\
        A_2'C_2 + b_{0,3,0}=0, \\
        A_2A_6' + A_6A_2' + A_2'C_3 +
        b_{0,2,1}=0, \\
        A_2'C_4 + b_{0,2,0}=0, \\
        A_3'C_2 +
        A_6'D_2 + b_{0,1,2}=0, \\
        A_3A_6' +
        A_6A_3' + A_3'C_3 + A_6'D_3 + b_{0,0,3}=0, \\
        A_3'C_4 + A_6'D_4 + b_{0,0,2}=0, \\
        A_3^2 + A_3A_3' + {A_3'}^2(\xi^2 +\xi) + a_{0,0,3}=0,\\
        A_3A_5' + A_5A_3' + A_3'D_1 + c_{0,0,3}=0, \\
        A_3A_6' + A_6A_3' + A_3'D_2 + b_{0,0,3}=0, \\
        A_3'D_3 + a_{0,0,3} + A_3'D_4 + a_{0,0,2}=0, \\
        A_1A_3' + A_3A_1' + A_5^2 + A_5A_5'
        + A_3'B_1 + {A_5'}^2(\xi^2 + \xi) +
        A_5'D_1 + D_1^2=0, \\ 
        A_3A_4' +
        A_4A_3' + A_5A_6' + A_6A_5' +
        A_3'B_2 + A_3'C_1 + A_51'D_2 +
        A_6'D_1=0,\\
        A_3A_5' + A_5A_3' +
        A_3'B_3 + A_3'D_1 + A_5'D_3=0,\\
        A_3'B_4 + A_5'D_4 + c_{0,0,2}=0,\\
        A_2A_3' +
        A_3A_2' + A_6^2 + A_6A_6' + A_3'C_2
        + {A_6'}^2(\xi^2 + \xi) + A_6'D_2 +
        D_2^2=0,\\
        A_3A_6' + A_6A_3'=0,\\
        A_3'C_3 + A_3'D_2 + A_6'D_3=0,\\
        A_3'C_4 + A_6'D_4 + b_{0,0,2}=0,\\
        A_3'D_3 +
        D_3^2=0,\\
        \end{eqnarray*}
        \begin{eqnarray*}
        A_3'D_4 + a_{0,0,2}=0,\\
        D_4^2 + a_{0,0,1}=0,\\
        A_1A_5' + A_5A_1' + A_1'D_1 + A_5'B_1 +
        a_{3,0,0}=0,\\
        A_1A_6' + A_4A_5' + A_5A_4'
        + A_6A_1' + A_1'D_2 + A_4'D_1 +
        A_5'B_2 + A_5'C_1 + A_6'B_1 +
        a_{2,1,0}=0,\\
        A_1'D_3 +A_5'B_3 +
        a_{2,0,1}=0,\\
        A_1'D_4 + A_5'B_4 +a_{2,0,0}=0,\\
        A_2A_5' + A_4A_6' + A_5A_2' +
        A_6A_4' + A_2'D_1 + A_4'D_2 +
        A_5'C_2 + A_6'B_2 + A_6'C_1 +
        a_{1,2,0}=0,\\
        A_4'D_3 + A_5'C_3 +
        A_6'B_3 + a_{1,1,1}=0,\\
        A_4D_4 +
        A_5'C_4 + A_6'B_4 + a_{1,1,0}=0,\\
        A_3A_5' +
        A_5A_3' + A_3'D_1 + a_{1,0,2}=0, \\
        a_{1,0,1}=0, \\
        a_{1,0,0} + c_{0,0,1}=0, \\
        A_2A_6' +A_6A_2' +
        A_2'D_2 + A_6'C_2 + a_{0,3,0}=0, \\
        A_2'D_3 +
        A_6'C_3 + a_{0,2,1}=0, \\
        A_2'D_4 + A_6'C_4 +
        a_{0,2,0}=0, \\
        A_3A_6' + A_6A_3' + A_3'D_2
        + a_{0,1,2}=0, \\
       A_3'D_3 + a_{0,0,3}=0, \\
        A_3'D_4 + a_{0,0,2}=0, \\
        A_1^2 + A_1A_1'+ {A_1'}^2(\xi^2 +\xi) + A_1'B_1
        + B_1^2 + c_{3,0,0}=0, \\
        A_1A_4' + A_4A_1' +
        A_1'B_2 + A_1'C_1 + A_4'B_1 + b_{3,0,0} +
        c_{2,1,0}=0, \\
        A_1A_5' + A_5A_1' + A_1'B_3 +
        A_1'D_1 + A_5'B_1 + a_{3,0,0} + c_{2,0,1}=0, \\
        A_1A_2' + A_2A_1' +
        A_4^2 + A_4A_4' + A_1'C_2 + A_2'B_1
        +\\ {A_4'}^2(\xi^2 + \xi) + A_4'B_2 +
        A_4'C_1 + B_2^2 + C_1^2 + b_{2,1,0} +
        c_{1,2,0}=0, \\
        A_1A_6' + A_4A_5' +
        A_5A_4' + A_6A_1' + A_1'C_3 +
        A_1'D_2 + A_4'B_3 +\\ A_4'D_1 +
       A_5'B_2 + A_5'C_1 + A_6'B_1 +
        a_{2,1,0} + b_{2,0,1} + c_{1,1,1}=0, \\
        A_1'C_4 +
        A_4'B_4 + b_{2,0,0} + c_{1,1,0}=0, \\
        A_1A_3' +
        A_3A_1' + A_5^2 + A_5A_5' + A_1'D_3
        + A_3'B_1 +\\ {A_5'}^2(\xi^2 + \xi) +
        A_5'B_3 + A_5'D_1 + B_3^2 + D_1^2 +
        a_{2,0,1} + c_{1,0,2}=0, \\
        A_1'D_4 + A_5'B_4 +
        a_{2,0,0} + c_{1,0,1}=0, \\
        A_2A_4' +
        A_4A_2' + A_2'B_2 + A_2'C_1 + A_4'C_2 +
        b_{1,2,0} + c_{0,3,0}=0, \\
        A_2A_5' + A_4A_6' +
        A_5A_2' + A_6A_4' + A_2'B_3 +
        A_2'D_1 + A_4'C_3 + A_4'D_2 +\\
        A_5'C_2 + A_6'B_2 + A_6'C_1 +
        a_{1,2,0} + b_{1,1,1} +c_{0,2,1}=0, \\
        A_2'B_4 +
        A_4'C_4 + b_{1,1,0} + c_{0,2,0}=0, \\
        A_3A_4' +
        A_4A_3' + A_5A_6' + A_6A_5' +
        A_3'B_2 + A_3'C_1 +\\ A_4'D_3 +
        A_5'C_3 + A_5'D_2 + A_6'B_3 +
        A_6'D_1 + a_{1,1,1} + b_{1,0,2} + c_{0,1,2}=0, \\
        \end{eqnarray*}
        \begin{eqnarray*}
        A_4'D_4 + A_5'C_4 + A_6'B_4 + a_{1,1,0} +
        b_{1,0,1} + c_{0,1,1}=0, \\
        A_3A_5'
        + A_5A_3' + A_3'B_3 + A_3'D_1 + A_5'D_3 +
        a_{1,0,2} + c_{0,0,3}=0, \\
        A_3'B_4 + A_5'D_4 +
        a_{1,0,1} + c_{0,0,2}=0, \\
        A_2^2 +
        A_2A_2' + {A_2'}^2(\xi^2 + \xi) + A_2'C_2 + C_2^2
        + b_{0,3,0}=0, \\
        A_2A_6' + A_6A_2' + A_2'C_3 +
        A_2'D_2 + A_6'C_2 + a_{0,3,0} + b_{0,2,1}=0, \\
        A_2A_3' + A_3A_2' +
        A_6^2 + A_6A_6' + A_2'D_3 + A_3'C_2
        +\\ {A_6'}^2(\xi^2 + \xi) + A_6'C_3 +
        A_6'D_2 + C_3^2 + D_2^2 + a_{0,2,1} +
        b_{0,1,2}=0, \\
         A_2'D_4 + A_6'C_4 + a_{0,2,0} +
        b_{0,1,1}=0, \\
        A_3A_6' + A_6A_3' +
        A_3'C_3 + A_3'D_2 + A_6'D_3 + a_{0,1,2} +
        b_{0,0,3}=0, \\
        A_3'C_4 + A_6'D_4 + a_{0,1,1} +
        b_{0,0,2}=0, \\
        A_3^2 + A_3A_3' +
       {A_3'}^2(\xi^2 + \xi) + A_3'D_3 + D_3^2 + a_{0,0,3}=0. \\
        \end{eqnarray*}

\end{document}